\documentclass[a4paper,10pt]{article}
\usepackage[margin=2cm]{geometry}
\usepackage[l2tabu, orthodox]{nag}
\usepackage{
    amsfonts,
    amsmath,
    amssymb,
}

\usepackage[usenames,dvipsnames]{xcolor}

\usepackage{newpxtext,newpxmath}

\usepackage[
    backend=biber,
    bibencoding=utf-8,
    bibwarn=true,
    style=alphabetic,
    sorting=nyt,
    doi=true,
    isbn=true,
    url=false,
    eprint=true,
    backref=true,
    maxnames=100,
]{biblatex}

\AtEveryBibitem{%
    \clearlist{language}%
}
\definecolor{blue-navy}{HTML}{002867}
\usepackage[
    colorlinks=true,
    linkcolor=blue-navy,
    anchorcolor=blue-navy,
    citecolor=blue-navy,
    urlcolor=blue-navy,
    plainpages=false,
    bookmarks=true,
    bookmarksopen=false,
]{hyperref}

\newcommand{\bbk}{{\mathbb{k}}}

\newcommand{\bbC}{{\mathbb{C}}}

\newcommand{\bbF}{{\mathbb{F}}}

\newcommand{\bbK}{{\mathbb{K}}}

\newcommand{\bbO}{{\mathbb{O}}}

\newcommand{\bbQ}{{\mathbb{Q}}}
\newcommand{\bbR}{{\mathbb{R}}}

\newcommand{\bbZ}{{\mathbb{Z}}}

\newcommand{\cH}{{\mathcal{H}}}

\newcommand{\cL}{{\mathcal{L}}}

\newcommand{\cN}{{\mathcal{N}}}

\newcommand{\cR}{{\mathcal{R}}}

\newcommand{\GL}{\operatorname{GL}}

\newcommand{\im}{\operatorname{im}}

\newcommand{\id}{\operatorname{id}}
\newcommand{\Hom}{\operatorname{Hom}}
\newcommand{\End}{\operatorname{End}}

\newcommand{\rank}{\operatorname{rank}}

\newcommand{\Mat}{\operatorname{Mat}}

\newcommand{\Sym}{\operatorname{Sym}}

\newcommand{\ch}{\operatorname{ch}}

\newcommand{\Std}{\operatorname{Std}}

\newcommand{\xsurjto}[1]{\xtwoheadrightarrow{#1}}

\newcommand{\injto}{\hookrightarrow}
\newcommand{\surjto}{\twoheadrightarrow}

\newcommand{\cat}[1]{\operatorname{\mathsf{#1}}}

\newcommand{\innprod}[1]{\langle #1 \rangle}

\newcommand{\abs}[1]{\lvert #1 \rvert}

\newcommand{\set}[1]{\left\{ #1 \right\}}

\usepackage[
    capitalise,
    nameinlink, 
]{cleveref}

\usepackage{blkarray}   
\usepackage{enumitem}   
\usepackage{cancel}     
\usepackage{extpfeil}   

\setlist[enumerate]{parsep=0pt}
\setlist[itemize]{parsep=0pt}

\newcommand{\defn}[1]{{\color{Brown}\emph{#1}}} 

\usepackage{tikz,tikz-cd}

\def\bscounit{\tikz[baseline=.1ex]{
    \draw[-] (0, -0.5ex) -- (0, 1.5ex) node[circle,fill=black,inner sep=0pt,minimum size=0.5ex] {};
}}
\def\bsunit{\tikz[baseline=.1ex]{
    \draw[-] (0, 1.5ex) -- (0, -0.5ex) node[circle,fill=black,inner sep=0pt,minimum size=0.5ex] {};
}}
\def\bsmult{\tikz[baseline=.1ex]{
    \draw[-] (1ex, -0.5ex) -- (0, 0.5ex) -- (0, 1.5ex);
    \draw[-] (-1ex, -0.5ex) -- (0, 0.5ex); 
}}

\newcommand{\LLup}{\operatorname{\mathsf{LL}}_{\mathrm{up}}}
\newcommand{\LLdown}{\operatorname{\mathsf{LL}}_{\mathrm{down}}}
\newcommand{\TLLup}{\operatorname{\mathsf{TLL}}_{\mathrm{up}}}
\newcommand{\TLLdown}{\operatorname{\mathsf{TLL}}_{\mathrm{down}}}

\newtheorem{theorem}{Theorem}[section]
\crefname{theorem}{Theorem}{Theorems}

\crefname{corollary}{Corollary}{Corollaries}
\newtheorem{lemma}[theorem]{Lemma}
\crefname{lemma}{Lemma}{Lemmas}
\newtheorem{remark}[theorem]{Remark}
\crefname{remark}{Remark}{Remarks}
\newtheorem{definition}[theorem]{Definition}
\crefname{definition}{Definition}{Definitions}

\usepackage{todonotes}

\title{Calculating the $p$-canonical basis of Hecke algebras}
\author{Joel Gibson, Lars Thorge Jensen, and Geordie Williamson \\ \\ \textit{\normalsize Dedicated to the memory of Jim Humphreys}}

\begin{document}

    \maketitle


\begin{abstract}
    We describe an algorithm for computing the $p$-canonical basis of the Hecke algebra, or one of its antispherical modules.
    The algorithm does not operate in the Hecke category directly, but rather uses a faithful embedding of the Hecke category inside a semisimple category to build a ``model'' for indecomposable objects and bases of their morphism spaces.
    Inside this semisimple category, objects are sequences of Coxeter group elements, and morphisms are (sparse) matrices over a fraction field, making it quite amenable to computations.
    This strategy works for the full Hecke category over any base field, but in the antispherical case we must instead work over $\bbZ_{(p)}$ and use an idempotent lifting argument to deduce the result for a field of characteristic $p > 0$.
    We also describe a less sophisticated algorithm which is much more suited to the case of finite groups.
    We provide complete implementations of both algorithms in the MAGMA computer algebra system.
\end{abstract}

\section{Introduction}

Let $(W, S)$ be a Coxeter system, which defines a Hecke algebra $H$, a $\bbZ[v^{\pm 1}]$-algebra with a basis indexed by $W$.
The Hecke algebra has a \textit{standard basis} which arises from viewing it as a deformation of the group ring $\bbZ[W]$.
Kazhdan and Lusztig \cite{kazhdanRepresentationsCoxeterGroups1979} defined \textit{canonical basis} of the associated Hecke algebra, and for each pair $y, x \in W$ the \textit{Kazhdan-Lusztig polynomial} $h_{y, x} \in \bbZ[v]$ as the entries appearing in the change of basis from the canonical to the standard basis.
A choice of subset $I \subseteq S$ defines a right module over the Hecke algebra with a basis indexed by cosets $W_I \backslash W$, which we call the \textit{antispherical module}.
This module also admits a canonical basis, and analagously for each $y, x \in W_I \backslash W$ there is an \textit{antispherical Kazhdan-Lusztig polynomial} $n_{y, x}^I \in \bbZ[v]$.

Since the initial work of Kazhdan and Lusztig \cite{kazhdanRepresentationsCoxeterGroups1979}, it has been realised \cite{kazhdanSchubertVarietiesPoincare1980,springerQuelquesApplicationsCohomologie1982} that the Hecke algebra admits a categorification, now referred to as the \textit{Hecke category}.
There is an isomorphism from the split Grothendieck group of this categorification to the Hecke algebra.
This isomorphism sends the classes of indecomposable objects to the canonical basis, implying many nice properties about the canonical basis (for example the positivity of the coefficients of the Kazhdan-Lusztig polynomials $h_{y, x}$).
The Hecke category may also be defined over a field of characteristic $p > 0$.
In this case the classes of indecomposable objects are sent to a different basis of the Hecke algebra, which is called the \textit{$p$-canonical basis}.
(We remark that the structure of the Hecke category in positive characteristic really depends on a choice of realisation for the Coxeter system $(W, S)$, a choice which does not matter in characteristic zero).
The polynomials ${^p} h_{y, x}$ which express the $p$-canonical basis in the standard basis are of great interest in modular representation theory.
For example, $p$-canonical basis elements appear in character formulas \cite{richeTiltingModulesPcanonical2017,richeSmithTreumannTheoryLinkage2020} for simple and tilting modules for reductive algebraic groups over a field of characteristic $p$.
They also control decomposition numbers for higher cyclotomic quotients of symmetric groups \cite{bowmanPathIsomorphismsQuiver2021}.

The Kazhdan-Lusztig polynomials $h_{y, x}$ and their antispherical analogues $n_{y, x}^I$ have an algebraic characterisation, which gives a straightforward inductive algorithm for their calculation.\footnote{
	The main obstacle for the computation of Kazhdan-Lusztig polynomials is the amount of RAM required: keeping track of all intermediate results (which seems somewhat necessary) requires $O(|W|^2)$ storage, which is super-factorial in the rank of a finite Coxeter group.
	See \cite{voganCharacterTableE82007} for the computational obstacles encountered during the calculation of similar Kazhdan-Lusztig-Vogan polynomials in type $E_8$.
}
In contrast, there is no known algebraic characterisation of the $p$-canonical polynomials, and it is unlikely that such a characterisation exists; for example, the results of \cite{williamsonSchubertCalculusTorsion2016} suggest that any such algorithm would have to involve number-theoretic considerations (for example, the divisibility of Fibonacci numbers by a fixed prime).
We give a practical (but slow) algorithm for the calculation of ${^p} h_{y, x}$, relying on the technique of localisation developed in the work of Elias and Williamson.
In the anti-spherical case, the algorithm makes use of ideas from Arikhipov-Bezrukavnikov [AB], Riche-Williamson \cite{richeSmithTreumannTheoryLinkage2020}, and Libedinsky-Williamson \cite{libedinskyAntisphericalCategory2017}.

An earlier (slightly different) incarnation of this algorithm was used to produce the computations leading to the Lusztig-Williamson ``billiards conjecture'' \cite{lusztigBilliardsTiltingCharacters2018}.
A version of the algorithm presented here was used by Jensen to generate data leading to his revised form of the billiards conjecture \cite{jensenCorrectionLusztigWilliamsonBilliards2021}.
The simpler algorithm presented at the end of this paper was used to generate the tables in the appendix of \cite{eliasCategoricalDiagonalizationCells2021}.

The algorithms are programmed using the computer algebra system MAGMA \cite{bosmaMagmaAlgebraSystem1997}.
The source code for the MAGMA package containing both algorithms is available online at \url{https://github.com/joelgibson/ASLoc}.

Organisation of the paper is as follows.
In \cref{section:background} we fix notation and briefly recall the definitions of Hecke algebra, realisation, root system, Hecke category, and antispherical category.
In \cref{section:localisation} we review localisation for the Hecke and antispherical categories, our key tool for computations.
Finally, the main algorithm is stated in \cref{section:algorithm}, and a simpler algorithm which works particularly well for finite groups of small rank is stated in \cref{section:simpler-algorithm}.

\subsection*{Dedication}

It is a great pleasure to dedicate this paper to the memory of Jim Humphreys.
His books on Coxeter groups and Lie algebras have become standard references, and were an early inspiration to us and many others.
We are deeply grateful for his contributions to our field.

\tableofcontents

\section{Background and notation} \label{section:background}

In this section we give an overview of the notation and conventions used throughout.
We direct the reader to \cite{eliasIntroductionSoergelBimodules2020} for a more detailed account of the diagrammatic Hecke category; our notation has been chosen to agree with this reference as much as possible.

\subsection{Overview of notation}

\begin{center}
    \begin{tabular}{rl}
        $(W, S)$                                                                                       & Coxeter system \\
        $\cR \colon W \to 2^S$                                                                         & Right descent set $\cR(w) = \{s \in S \mid l(ws) < l(w)\}$ \\
        $(W_I, I)$                                                                                     & Parabolic sub-Coxeter system defined by $I \subseteq S$ \\
        ${^I W}$                                                                                       & Minimal right coset representatives for $W_I \backslash W$ \\
        $p$                                                                                            & A prime or zero \\
        $\mathbb{Z}_{(p)}$                                                                             & $\bbZ$ localised at $(p)$ \\
        $(\bbO, \bbK, k)$                                                                              & A $p$-modular system, usually $(\bbZ_{(p)}, \bbQ, \bbF_p)$ \\
        $\Delta = \{\alpha_s\}_{s \in S} \subseteq X$                                                  & Set of simple roots in the free module $X$ \\
        $\Delta^\vee = \{\alpha_s^\vee\}_{s \in S} \subseteq X^\vee$                                   & Set of simple coroots in the free module $X^\vee$ \\
        $(\Delta^\vee \subseteq X^\vee, \Delta \subseteq X)$                                           & Realisation of $(W, S)$ over $\bbZ$, assumed to be free \\
        $\Phi \subseteq X$                                                                             & Root system for $W$: $\Phi = \{w(\alpha_s) \mid w \in W, s \in S\}$ \\
        $\Phi_I$                                                                                       & Root system for $W_I$: $\Phi_I = \{w(\alpha_s) \mid w \in W_I, s \in I\}$ \\
        $R = \Sym^\bullet(X)$                                                                          & Polynomial ring with $X$ in degree 2 \\
        $Q = R[\beta^{-1} \mid \beta \in \Phi]$                                                      & Localisation of $R$ with all roots inverted \\
        $R_I = R / \innprod{\alpha_s \mid s \in I}$                                                    & The $I$-quotient ring \\
        $R_{(I)} \subseteq Q$                                                                          & Subring of $Q$ where the $\Phi_I$ are \textit{not} inverted \\
        $Q_I = R_I[\beta^{-1} \mid \beta \in \Phi \setminus \Phi_I]$                                 & Localisation of $R_I$ with $\Phi \setminus \Phi_I$ inverted
    \end{tabular}
\end{center}

\subsection{Hecke algebras and the antispherical module} \label{section:hecke-algebras}

Let $(W, S)$ be a Coxeter system, with $[m_{st}]$ the associated Coxeter matrix, and let $\cL = \bbZ[v^{\pm 1}]$ denote the ring of Laurent polynomials with integer coefficients.
The \defn{Hecke algebra} $H = H(W, S)$ is the unital $\cL$-algebra generated by the elements $\{\delta_s \mid s \in S\}$, subject to the two relations
\begin{enumerate}
    \item (\textit{Quadratic relation}): $\delta_s^2 = (v^{-1} - v)\delta_s + 1$ for all $s \in S$, and
    \item (\textit{Braid relation}): $\underbrace{\delta_s \delta_t \delta_s \cdots}_{m_{st} \text{ terms}} = \underbrace{\delta_t \delta_s \delta_t \cdots}_{m_{st} \text{ terms}}$ for all $s, t \in S$ with $m_{st} < \infty$.
\end{enumerate}
Let $(s_1, \ldots, s_n)$ be a reduced expression for some $x \in W$ in the generators $S$, then define the element $\delta_x = \delta_{s_1} \cdots \delta_{s_n}$.
The element $\delta_x$ is independent of the choice of reduced expression, and the set $\{\delta_x \mid x \in W\}$ forms a basis for the Hecke algebra, called the \defn{standard basis}.
Each standard basis element is invertible, since each $\delta_s$ is invertible by the quadratic relation.

The \defn{bar involution} on $H$ is the unique $\bbZ$-linear map of algebras $\sigma \colon H \to H$ satisfying $\sigma(\delta_x) = \delta_{x^{-1}}^{-1}$ and $\sigma(v) = v^{-1}$.
The \defn{Kazhdan-Lusztig basis} or \defn{canonical basis} is the unique set $\{b_x \mid x \in W\}$, where each element is fixed under $\sigma$ and further satisfies the Bruhat upper-triangularity and degree-bound condition
\begin{equation}
    b_x = \delta_x + \sum_{y < x} h_{y, x} \delta_y, \quad \text{where } h_{y, x} \in v \bbZ[v].
\end{equation}
The polynomials $h_{y, x}$ are called \defn{Kazhdan-Lusztig polynomials}.
We define $h_{x, x} = 1$ and $h_{y, x} = 0$ if $y \not \leq x$.

For any subset $I \subseteq S$ of the generators, the subgroup $W_I$ generated by $I$ is called the \defn{standard parabolic subgroup} associated to $I$.
Each right coset $W_I x \in W_I \backslash W$ admits a unique representative of minimum length; let ${^I W} \subseteq W$ be this set of \defn{minimal coset representatives}.
The set of minimal coset representatives may be described in several convenient ways (see \cref{section:root-systems} for the definition of a root system):
\begin{equation}
    {^I W}
    = \{w \in W \mid sw > w \text{ for all } s \in I\}
    = \{w \in W \mid w^{-1}(\beta) > 0 \text{ for all } \beta \in \Phi^+_I\}.
\end{equation}
The multiplication map $W_I \times {^I W} \to W$ is a bijection of sets, and furthermore is additive in length: $l(uv) = l(u) + l(v)$ for $u \in W_I$ and $v \in {^I W}$.
The pair $(W_I, I)$ is a Coxeter system in its own right and defines a Hecke algebra $H_I = H(W_I, I)$, which is a subalgebra of $H$ because of the length-additivity.

Due to the quadratic relation, the map $\varphi_I \colon H_I \to \cL$ defined on generators by $\delta_s \mapsto -v$ extends to an algebra homomorphism, making $\cL$ into an $(H_I, H_I)$-bimodule we denote by $\cL(-v)$.
The \defn{antispherical module} ${^I N}$ is the right $H$-module defined via tensor induction:
\begin{equation}
    {^I N} = \cL(-v) \otimes_{H_I} H.
\end{equation}
The set $\{1 \otimes \delta_x \mid x \in {^I W}\}$ forms a basis for ${^I N}$, called the \defn{standard basis}.
The bar involution $\sigma \colon H \to H$ induces a map of $\bbZ$-modules $\sigma_I \colon {^I N} \to {^I N}$ by setting $\sigma_I(p \otimes h) = \sigma(p) \otimes \sigma(h)$.
The \defn{canonical basis} of the antispherical module ${^I N}$ is the unique set $\{d_x \mid x \in {^I W}\} \subseteq {^I N}$ which is fixed pointwise under $\sigma_I$, and satisfies the degree bound condition
\begin{equation}
    d_x = 1 \otimes \delta_x + \sum_{\substack{y < x \\ y \in {^I W}}} n_{y, x} (1 \otimes \delta_y), \quad \text{where } n_{y, x} \in v \bbZ[v].
\end{equation}
The polynomials $n_{y, x}$ are called \defn{antispherical Kazhdan-Lusztig polynomials}.
We again define $n_{x, x} = 1$ for all $x \in {^I W}$, and $n_{y, x} = 0$ for $x, y \in {^I W}$ with $y \not\leq x$.

\subsection{Realisations of Coxeter systems}

While the Hecke algebra and canonical basis depend only on the Coxeter system $(W, S)$, the $p$-canonical basis and the Hecke category depend on the data of a \textit{realisation} of the Coxeter system, which comes with a choice of field: it is the characteristic of this field which gives the $p$-canonical basis its $p$.
Take for example the Cartan data $B_n$ and $C_n$ for $n \geq 3$, which define identical Coxeter systems and Hecke algebras, but different $p$-canonical bases for $p = 2$, as was shown in \cite{jensenCanonicalBasisHecke2016}.
We are primarily interested in the $p$-canonical bases for the finite and affine Weyl groups, whose root data are defined over the integers, and hence over any characteristic by base change.
Below we review some standard theory about realisations of root systems over the integers, taking generalised Cartan matrices as a starting point.

Let $A = [a_{ij}] \in \Mat_n(\bbZ)$ be a \defn{generalised Cartan matrix} (GCM), meaning that (a) $a_{ii} = 2$ for all $i$, (b) $a_{ij} \leq 0$ for all $i \neq j$, and (c) $a_{ij} = 0 \iff a_{ji} = 0$.
Each GCM has an associated Coxeter system $(W, S)$, which is determined by the products $a_{ij} a_{ji}$ of off-diagonal elements:
for $i \neq j$, let $m_{ij} = 2, 3, 4, 6$, or $\infty$ depending on whether $a_{ij}a_{ji}$ is equal to $0, 1, 2, 3$, or $\geq 4$ respectively.
Let $(W, S)$ be the Coxeter system with Coxeter matrix $[m_{ij}]$.

A \defn{realisation} of $A$ over $\bbZ$ is a pair $X^\vee, X$ of free $\bbZ$-modules of finite rank, together with a perfect pairing $\innprod{-, -} \colon X^\vee \times X \to \bbZ$ and sets $\Delta^\vee = \set{\alpha_s^\vee \mid s \in S} \subseteq X^\vee$ and $\Delta = \set{\alpha_s \mid s \in S} \subseteq X$ such that $\innprod{\alpha_s^\vee, \alpha_t} = a_{st}$.
The $\alpha_s^\vee$ are called the \defn{simple coroots}, and the $\alpha_s$ are called the \defn{simple roots}.
For any realisation, the maps $r_s^\vee \in \GL_\bbZ(X^\vee)$ and $r_s \in \GL_\bbZ(X)$ defined by
\begin{align}
    r_s^\vee(\mu) &= \mu - \innprod{\mu, \alpha_s}\alpha_s^\vee
    &
    r_s(\lambda) &= \lambda - \innprod{\alpha_s^\vee, \lambda} \alpha_s
\end{align}
satisfy the braid relations for $(W, S)$, and hence give representations of $(W, S)$ on $X^\vee$ and $X$.
These representations are dual in the sense that $\innprod{w \nu, w \lambda} = \innprod{\nu, \lambda}$ for all $w \in W$, $\nu \in X^\vee$, and $\lambda \in X$ (so if dual bases are chosen for $X$ and $X^\vee$, the matrices associated to $w$ acting on each side are transpose).

A realisation is called \defn{cofree} if the simple coroots are linearly independent, and \defn{free} if the simple roots are linearly independent.
If a realisation is free or cofree, then the $W$-action on $X^\vee$ and $X$ is faithful.
If $A$ is of finite type any realisation is automatically both free and cofree, whereas if $A$ is of affine type then a realisation may be free without being cofree, or vice-versa.

There is a canonical way of producing a free realisation of any GCM $A$, called the \defn{adjoint} realisation\footnote{
    This terminology comes from the theory of split reductive groups: the \textit{adjoint} realisation as defined here is the root data associated to the reductive group of adjoint type for the generalised Cartan matrix $A$.
}.
Define $X = \bigoplus_{s \in S} \bbZ \alpha_s$ to be the free $\bbZ$-module with basis $\{\alpha_s\}_{s \in S}$ of simple roots, and set $X^\vee = \Hom(X, \bbZ)$, defining the simple coroots $\alpha_s^\vee$ as the unique solutions to the system of linear equations $\innprod{\alpha_s^\vee, \alpha_t} = a_{st}$.
Concretely, the simple roots $\alpha_s$ may be taken to be coordinate column vectors, the simple coroots $\alpha_s^\vee$ to be the rows of the Cartan matrix $A$, and the perfect pairing $\innprod{-, -}$ to be the dot product.

A realisation satisfies \defn{Demazure surjectivity} if the two maps
\begin{align}
    \innprod{\alpha_s^\vee, -} &\colon X \to \bbZ
    &
    \innprod{-, \alpha_s} &\colon X^\vee \to \bbZ
\end{align}
are surjective for all $s \in S$.
If a realisation satisfies Demazure surjectivity over $\bbZ$, then after base change to any field it also satisfies Demazure surjectivity.
Furthermore, after base change to a field of characteristic not 2, every realisation satisfies Demazure surjectivity because in particular $\innprod{\alpha_s^\vee, \alpha_s} = 2$ is invertible.
Demazure surjectivity is necessary for having a well-behaved Hecke category, since the double leaves may fail to span the hom space between two Bott-Samelson bimodules \cite[\S~5.2]{eliasLocalizedCalculusHecke2020} (the double leaves remain linearly independent regardless of whether the realisation satisfies Demazure surjectivity).

Our algorithms use the adjoint realisation, which may fail to be Demazure-surjective.
For example in the adjoint realisation of type $A_1$ (with Cartan matrix $A = (2)$), the image of $\innprod{\alpha^\vee, -}$ is $2 \bbZ$.
We briefly justify why the algorithm remains correct despite not using a Demazure-surjective realisation.

Any realisation may be extended to one in which Demazure surjectivity holds.
A natural example of this in type $A_1$ is to use the root data associated to $\GL_2$, in which $X = X^\vee = \bbZ^2$, the perfect pairing is the dot product, and $\alpha_1 = \alpha_1^\vee = (1, -1)$.
In order to extend the adjoint realisation to one in which Demazure surjectivity holds, add extra basis elements $d$ to $X$ and $d^\vee$ to $X^\vee$, and extend the perfect pairing in such a way that $\innprod{\alpha_s^\vee, d} = 1$ for all $s \in S$.
We end up with a map $f \colon X \to X^e$ from the adjoint realisation into the extended realisation which preserves simple roots (along with a transpose map $g \colon (X^e)^\vee \to X^\vee$ preserving simple coroots).
By functorality we get a map $R \to R^e$, where $R = \Sym(X)$ is the usual polynomial ring in the definition of the Hecke category.

In our Hecke category defined over $R$, we end up with light leaves with coefficients in $R$, which are independent over both $R$ and $R^e$, but may fail to span correctly for the Hecke category defined over $R$.
However this is not a problem: we may imagine our light leaves as sitting inside the Hecke category defined over $R^e$ where they do span correctly (and just happen to have coefficients in $R$).
It is due to the explicit light leaves construction that we can do this: some other way to find a basis which depended only on the data of $R$ might find some extra elements.

\subsection{Root systems and parabolic subgroups} \label{section:root-systems}

Let $(X^\vee, X, \Delta^\vee, \Delta)$ be a realisation over the integers, or an ordered characteristic zero ring such as $\bbR$ or $\bbZ_{(p)}$.
The elements in the orbit $\Phi = \set{w(\alpha_s) \mid w \in W, s \in S} \subseteq X$ of the simple roots under the action of $W$ are called \defn{roots}.
Provided that the realisation is free\footnote{
    A weaker condition suffices: the simple roots should be \defn{positively independent}, meaning that if $\sum_{s \in S}\lambda_s \alpha_s = 0$ for $\lambda_s \geq 0$, then all $\lambda_s = 0$.
    This is equivalent to the existence of an $f \in X^\vee$ such that $\innprod{f, \alpha_s} > 0$ for all $s \in S$.
},
these roots may be partitioned into two disjoint sets of positive roots $\Phi^+ = \Phi \cap \sum_{s \subseteq S} \bbZ_{\geq 0} \alpha_s$ and negative roots $-\Phi^+$, in which case we say that the $\Phi$ is a \defn{root system} for the realisation.
A root $\beta \in \Phi$ is called a \defn{positive root} if $\beta \in \Phi^+$, which we also write as $\beta > 0$; negative roots are defined similarly.

A standard parabolic subgroup $(W, I)$ defines a subset $\Delta_I$ of the simple roots, and a subset $\Phi_I \subseteq \Phi$ of the roots called the \defn{parabolic roots}:
\begin{equation}
    \Phi_I = \{w(\alpha_s) \mid w \in W_I, s \in I\}.
\end{equation}
The positive roots of $\Phi_I$ are the positive roots of $\Phi$ belonging to $\Phi_I$, so $\Phi_I^+ = \Phi^+ \cap \Phi_I$, and likewise for the negative roots $\Phi_I^-$.
The subset $\Phi_I$ may be written as the intersection of the $\bbZ$-span of $\Delta_I$ with the roots $\Phi$, so the root system $\Phi_I$ is cut out by a linear subspace in $X$.

The reader may recall from standard texts such as \cite{humphreysReflectionGroupsCoxeter1990} that positive and negative roots can be used to detect whether right multiplication by a generator increases or decreases the length of an element in a Coxeter system $(W, S)$: the precise statement is that
\begin{equation}
    l(ws) = \begin{cases}
        l(ws) = l(w) - 1 & \text{if } w(\alpha_s) < 0, \\
        l(ws) = l(w) + 1 & \text{if } w(\alpha_s) > 0,
    \end{cases}
\end{equation}
for all $w \in W$ and $s \in S$.
When working with the set
\[{^I W} = \set{w \in W \mid sw > w \text{ for all } s \in I}\]
of minimal right coset representatives for $W_I \backslash W$, there are three cases for right multiplication by $s \in S$: for a representative $x \in {^I W}$, the product $xs$ may be again in ${^I W}$, of length $l(x) \pm 1$, or $xs$ may be in the same coset as $x$ and no longer a minimal representative.
This trichotomy of cases is described in the following lemma.

\begin{lemma} \label{result:parabolic-property}
    \defn{(Parabolic property)} With notation as above, for $x \in {^I W}$ and $s \in S$ we have:
    \begin{equation}
        \begin{cases}
            xs \in {^I W} \text{ and } l(xs) = l(x) - 1 & \text{if } x(\alpha_s) \notin \Phi_I \text{ and } x(\alpha_s) < 0, \\
            xs \in {^I W} \text{ and } l(xs) = l(x) + 1 & \text{if } x(\alpha_s) \notin \Phi_I \text{ and } x(\alpha_s) > 0, \\
            xs \notin {^I W} & \text{if } x(\alpha_s) \in \Phi_I \text{ (if and only if } x(\alpha_s) \in \Delta_I \text{)}.
        \end{cases}
    \end{equation}
    Furthermore, in the last case we have $xs = tx$ for the unique $t \in I$ with $x(\alpha_s) = \alpha_t$.
\end{lemma}

A proof of the parabolic property for root systems over $\bbR$ may be found in 2.3 of \cite{libedinskyAntisphericalCategory2017}, the argument carries over verbatim for free realisations over $\bbZ$.

\subsection{The Hecke category}

Let $(W, S)$ be a Coxeter system, with a fixed choice $(X, X^\vee)$ of a realisation over $\bbZ$.
Let $R_\bbZ = \Sym_\bbZ^\bullet(X)$ be the polynomial algebra, graded with $X$ in degree 2; the action of $W$ on $X$ extends to $R_\bbZ$ by functorality of $\Sym$.
The \defn{diagrammatic Bott-Samelson category} $\cH_{BS}$ is a strict monoidal category enriched over $R_\bbZ$-modules, the definition of which is given in \cite[\S 10.2.4]{eliasIntroductionSoergelBimodules2020}.
We assume from this point on that the reader is familiar with the standard notation used for Soergel bimodules and the diagrammatic Hecke category which constitute parts I and II of \cite{eliasIntroductionSoergelBimodules2020}.

For any commutative ring $\bbk$, base change of the realisation $(X, X^\vee)$ to $\bbk$ commutes with base change of the category $\cH_{BS}$ to $\bbk$ given by $\cH_{BS} \otimes \bbk$, hence we get a diagrammatic Bott-Samelson category defined over $\bbk$.
The \defn{Hecke category} $\cH_\bbk$ is the additive, Karoubi envelope (idempotent completion) $Kar(\cH_{BS} \otimes \bbk)$.

The Hecke category $\cH_\bbk$ behaves best when $\bbk$ is a complete local ring, for instance a field or the $p$-adic integers $\bbZ_p$. In this case:
\begin{enumerate}
    \item The category $\cH_\bbk$ is Krull-Schmidt, in particular an object is indecomposable if and only if its endomorphism ring is local, and every object decomposes uniquely into a sum of indecomposables.
    \item The indecomposable objects $\set{{^\bbk B}_w \mid w \in W}$ are, up to isomorphism and grading shift, indexed by the elements of $W$.
    \item The \defn{character map} $\ch \colon [\cH_\bbk] \to H$ taking the class $[{^\bbk B}_s]$ to $b_s$ gives an isomorphism of algebras from the split Grothendieck group of $\cH_\bbk$ to the Hecke algebra.
        Hence we say that the Hecke category \defn{categorifies} the Hecke algebra.
\end{enumerate}
For $\bbk$ a complete local ring of characteristic $p$, we define the \defn{$p$-canonical basis} element ${^ p}b_w = \ch [{^\bbk B}_w]$.
It was shown in \cite{eliasHodgeTheorySoergel2014} that for $\bbk$ a field of characteristic zero we have ${^0 b_w} = b_w$, i.e. the $0$-canonical basis is identical to the canonical basis as defined in \cref{section:hecke-algebras}.

The character map is defined by the condition above, but is more transparently calculated by other means.
Let $\cH_\bbk^{\not <w}$ be the quotient of the category $\cH_\bbk$ by morphisms factoring through the additive subcategory generated by the objects $\set{{^\bbk B}_x \mid x < w}$ for $x$ lower in the Bruhat order than $w$.
As was shown in \cite[\S 6.3]{eliasSoergelCalculus2016}, the character map is given by
\begin{equation} \label{equation:character-map}
    \ch [B] = \sum_{w \in W} \underline{\rank} \Hom_{\not < w}^\bullet(B, {^\bbk B}_w) b_w,
\end{equation}
where by $\underline{\rank}$ we mean the graded rank as free $R_\bbk$ modules\footnote{
    If a $\bbZ$-graded vector space $V$ consists of two dimensions in degree $-1$ and three dimensions in degree $5$, then its graded rank is $\underline{\rank}(V) = 2v^{-1} + 3v^5$.
}.

\begin{remark}
    While the formation of the Bott-Samelson category $\cH_{BS}$ commutes with the base change of the realisation, the formation of the Hecke category $\cH$ does not, because additional idempotents might be present in $\cH_\bbk$ which are not present in $\cH_\bbZ \otimes \bbk$.
    In formulas, we have $\cH_\bbK = Kar(\cH_{BS} \otimes \bbk) = Kar(Kar(\cH_{BS}) \otimes \bbk)$.
    We also remark that the Hecke category $\cH$ over $\bbZ$ has no Krull-Schmidt property, and its indecomposable objects behave in quite interesting ways, which are not understood.

    We remark also that the endomorphism rings of ${^\bbk B}_x$ for $\bbk = \bbF_p$ are absolutely local, that is they remain local after any base change.
    Consequently, formation of the Hecke category does commute with base change from $\bbF_p$ to any extension, and so we may refer to the $p$-canonical basis rather than the $\bbk$-canonical basis.
\end{remark}

We spend the remainder of this section discussing light leaves, double leaves, and truncated leaves, which are our method at getting a grasp on the morphism spaces in $\cH_\bbk$.
We highly recommend the reader re-familiarise themselves with description of light leaves in \cite[\S~10.4]{eliasIntroductionSoergelBimodules2020} before moving on.

The objects of $\cH_{BS}$ are words in the Coxeter generators $S$.
Given a reduced expression $\underline{w}$ (in other words, an object of $\cH_{BS}$), an element $x \leq w$, and a fixed choice $\underline{x}$ of reduced expression for $x$, there are explicit sets of \defn{light leaves down} and \defn{light leaves up}
\begin{equation}
    \LLdown(\underline{w}, \underline{x}) \subseteq \Hom_{\cH_{BS}}(\underline{w}, \underline{x}), \quad
    \LLup(\underline{x}, \underline{w}) \subseteq \Hom_{\cH_{BS}}(\underline{x}, \underline{w}),
\end{equation}
each of which is in bijection with the set $\set{\underline{e} \subseteq \underline{w} \mid w^e = x}$ of subexpressions of $\underline{w}$ evaluating to $x$.
Given two reduced expressions $\underline{v}$ and $\underline{w}$, then after having fixed a choice of reduced expression for all $x \in W$ with $x \leq v$ and $x \leq w$, the morphisms from $\underline{v}$ to $\underline{w}$ have a \defn{double leaves basis}
\begin{equation}
    \bigsqcup_{\substack{x \leq v \\ x \leq w}} \LLup(\underline{x}, \underline{w}) \circ \LLdown(\underline{v}, \underline{x}) \subseteq \Hom_{\cH_{BS}}(\underline{v}, \underline{w}).
\end{equation}
The double leaves basis is primarily of theoretical interest: our algorithms will never actually compose two light leaves in this way to get a double leaf, since all the information is already contained just in the down and up leaves.

The light leaves can be put together into the double leaves basis for $\cH_{BS}$, while something similar can be done for \textit{truncated} light leaves in $\cH_{\bbk}$.
The indecomposable object ${^\bbk B}_w$ appears uniquely as a direct summand in $B_{\underline{w}}$, where $\underline{w}$ is a reduced expression for $w$, hence we get uniquely defined \textit{projection} and \textit{inclusion} morphisms
\begin{equation}
    p_{\underline{w}} \colon B_{\underline{w}} \surjto B_w, \quad i_{\underline{w}} \colon B_w \injto B_{\underline{w}},
\end{equation}
such that $e_{\underline{w}} := i_{\underline{w}} \circ p_{\underline{w}} \in \End(B_{\underline{w}})$ is idempotent, and $p_{\underline{w}} \circ i_{\underline{w}}$ is the identity on ${^\bbk B}_w$.
(Part of this statement is tautological: since $\cH_\bbk$ is defined as an idempotent completion of $\cH_{BS} \otimes \bbk$, the object ${^\bbk B_w}$ would really be formally defined as the pair $(B_{\underline{w}}, e_{\underline{w}})$ up to isomorphism, with $p_{\underline{w}}$ and $i_{\underline{w}}$ the associated projection and inclusion).
From the sets $\LLdown(\underline{w}, -)$ and $\LLup(-, \underline{w})$ we get sets of \defn{truncated light leaves}:
\begin{align}
    \TLLdown({^\bbk B}_w, {^\bbk B}_x) \subseteq p_{\underline{x}} \circ \LLdown(\underline{w}, \underline{x}) \circ i_{\underline{w}},
    \\
    \TLLup({^\bbk B}_x, {^\bbk B}_w) \subseteq p_{\underline{w}} \circ \LLdown(\underline{x}, \underline{w}) \circ i_{\underline{x}},
\end{align}
by picking some pair of reduced expressions, precomposing and postcomposing with the relevant inclusions and projections, and finally selecting a linearly independent subset (since in general many of the light leaves will have become linearly dependent after \textit{truncating}, i.e. precomposing and postcomposing by inclusions/projections).

We emphasise to the reader that the truncated light leaves are highly noncanonical, in three obvious ways:
\begin{enumerate}
    \item Neither the light leaves nor double leaves are canonical in $\cH_{BS}$.
    \item The reduced expressions chosen for $\underline{x}$ and $\underline{w}$ will affect the set of truncated leaves.
    \item We do not specify which morphisms to remove to make $\TLLdown$ and $\TLLup$ linearly independent.
\end{enumerate}

Assuming that we have picked a reduced expression for each element of the Coxeter group, the truncated light leaves can be put together into a set of truncated double leaves
\begin{equation}
    \Hom_{\cH_\bbk}({^\bbk B_x}, {^\bbk B_z}) = \bigoplus_{\substack{y \leq x \\ y \leq z}} \TLLup(y, z) \circ \TLLdown(x, y).
\end{equation}

\begin{remark}
    In order to calculate a full set of truncated light leaves, the difficulty lies in calculating the idempotent $e_{\underline{w}}$, or in other words trying to split lower terms off the Bott-Samelson object $B_{\underline{w}}$.
    The Bott-Samelson object $B_{\underline{w}}$ looks the same in every Hecke category (it simply base changes from the corresponding object in $\cH_{BS}$), but the indecomposable objects ${^\bbk B}_w$ depend on the field $\bbk$, and lower objects (and their corresponding truncated leaves) need to be calculated before higher idempotents can be calculated.
\end{remark}

\subsection{The antispherical category}

The \defn{antispherical category} ${^I \cN_\bbk}$ is defined either by taking the quotient (as additive categories) of $\cH_{BS}$ by all morphisms factoring through an element starting with an element of $I$, and then taking the base change to $\bbk$ followed by the Karoubi completion, or alternatively by taking the quotient of $\cH_\bbk$ by the right $\otimes$-ideal of all objects $\set{B_x \mid x \notin {^I W}}$: these two definitions agree by Proposition 3.2 of \cite{libedinskyAntisphericalCategory2017}.

\begin{remark}
    The equivalence of these categories is hinted at on the level of Hecke algebras and antispherical modules, as the antispherical module ${^I N}$ is naturally isomorphic to the quotient of $H$ by the right ideal generated by $\set{b_x \mid x \notin {^I W}}$.
\end{remark}

The antispherical category is a right module over the Hecke category, in the sense of there being a monoidal action functor
${^ I \cN_\bbk} \otimes \cH_\bbk \to {^ I \cN_\bbk}$: this functor canonically arises by the monoidal structure on $\cH_\bbk$ and the definition of the antispherical category as a quotient.
(This is a categorification of the fact that a ring modulo a right ideal is a right module).
The antispherical category ${^I \cN_\bbk}$ enjoys many properties analagous to that of the Hecke category: there is a similar \textit{antispherical} double-leaves basis that can be used to explicitly describe morphism spaces between Bott-Samelson objects, and if $\bbk$ is a field or complete local ring then the indecomposable objects are in bijection with minimal coset representatives ${^I W}$.
There is an analagous character map giving an isomorphism between the Grothendieck group of ${^I \cN_\bbk}$ and the antispherical module ${^I N}$, and for $\bbk$ a complete local ring of characteristic $p$, the image of the indecomposable objects is defined to be the \defn{$p$-canonical basis} of the antispherical module ${^I N}$.

\begin{remark}
    The antispherical category ${^I \cN}_\bbk$ has left and right actions of $R$, just as the Hecke category does.
    However the left action is simpler in the antispherical case, as multiplication by any $\set{\alpha_s \mid s \in I}$ on the left acts by zero.
    This can be seen diagrammatically: if a diagram of $\cH_{BS}$ has an $\alpha_s$-coloured strand appearing on the left for some $s \in I$, then the diagram becomes zero in the quotient.
    Since multiplication by $\alpha_s$ is appending an $\alpha_s$-coloured barbell on the left of the diagram, the result is zero in the quotient.
    So we may consider ${^I \cN}_\bbk$ to be an $(R/\innprod{\alpha_s \mid s \in I}, R)$-linear category.
\end{remark}


\section{Localisation} \label{section:localisation}

\defn{Localisation}, as we discuss it below, is the process of embedding a complicated category (the Hecke category, defined as as the additive Karoubi envelope of $\cH_{BS}$, itself a category presented by generators and relations) into a simple category (matrices with a twist).
Thus difficult questions become simple to answer: rather than ``do these two diagrams define the same morphism in the quotient?'' we ask ``are these two matrices equal?'', a question for which computer algebra software is well-equipped.
Localisation is our key tool to mechanically perform calculations in the Hecke category.

Some adjustments need to be made in the antispherical case, especially in positive characteristic where we need to use a $p$-modular system rather than just a field of characteristic $p$.
We encourage the reader have \cite[\S~1--2]{eliasLocalizedCalculusHecke2020} open beside them as they read our account of localisation.

\subsection{The standard category} \label{section:standard-category}

Before stating how localisation works we will define the target of localisation: a semisimple category we call the \textit{standard category}.
In what follows, $Q$ will be a ring admitting a natural map from $R$, for example in the case of the Hecke category $Q$ will be the field of fractions of $R$.

\begin{definition} \label{defn:q-groupoid}
    Let $W$ be a set, $Q$ a commutative ring, and $\varphi \colon W \to \End(Q)$ a family of ring endomorphisms.
    To this data we associate the \defn{$Q$-groupoid category} $\Omega(W, Q, \varphi)$ which is defined as follows:
    \begin{itemize}
        \item The objects are in bijection with $W$, labelled $\{r_x \mid x \in W\}$.
        \item The morphism spaces are left $Q$-modules, with $\End(r_x)$ free of rank $1$ spanned by $1_x$, and $\Hom(r_x, r_y) = 0$ for $x \neq y$.
        \item The morphism spaces are also right $Q$-modules, with the twisted action $1_x \cdot f = \varphi_x(f) 1_x$ for $f \in Q$, making each $\End(r_x)$ into a $(Q, Q)$-bimodule
    \end{itemize}

    If $W$ is a monoid and $\varphi$ a morphism of monoids, then setting $r_x \otimes r_y = r_{xy}$ and $f 1_x \otimes g 1_y = f \varphi_x(g) 1_{xy}$ (tensor product of bimodules over $Q$) defines a monoidal structure on $\Omega(W, Q, \varphi)$.
    If $W$ is a group, then every object is invertible under the monoidal structure.
\end{definition}

\begin{remark}
    The term \textit{groupoid} usually refers to a category in which every morphism is invertible, which the $Q$-groupoid is not.
    Rather, every object in the category is invertible under the monoidal product.
    These are sometimes called Picard categories in the literature, as the category of line bundles on a scheme shares this property.
\end{remark}

The category $\Omega(W, Q, \varphi)_{\oplus}$ is the \defn{additive envelope} of $\Omega(W, Q, \varphi)$, whose objects are formal sequences of objects of $\Omega(W, Q, \varphi)$, with morphisms being matrices whose entries are themselves morphisms in $\Omega(W, Q, \varphi)$.
(Thus, we make $\Omega(W, Q, \varphi)_\oplus$ into an additive category in the obvious way.)
If $W$ is a monoid then $\Omega(W, Q, \varphi)_{\oplus}$ also inherits a monoidal structure.

\begin{definition}
    The \textit{standard category} is $\Std_Q = \Omega(W, Q, \varphi)_{\oplus}$, where $W$ is the Coxeter group, $Q$ the field of fractions of $R$, and $\varphi$ the natural action of $W$ on $Q$.
\end{definition}

Since all morphisms in $\Omega(W, Q, \varphi)$ are either $0$ or scalar multiples of the identity, a morphism $\varphi \colon (x_1, \ldots, x_m) \to (y_1, \ldots, y_n)$ in the additive envelope $\Std_Q$ may be identified with an $n \times m$ matrix with values in $Q$ (we stress that the morphism is really the data of the domain sequence, codomain sequence, \textit{and} matrix, rather than just the matrix).
Composition in $\Std_Q$ becomes matrix multiplication, while the tensor product in $\Std_Q$ becomes a $\varphi$-twisted version of the tensor product of matrices.

In the $Q$-groupoid category there is only the zero homomorphism between $r_x$ and $r_y$ for $x \neq y$.
Therefore in $\Std_Q$ matrices will tend to be sparse, since a matrix entry is forced to be zero whenever its domain and codomain elements are distinct.

\textit{An example:} Let $W = \innprod{s, t \mid s^2 = t^2 = (st)^3 = 1}$ be the symmetric group on three letters, and $Q = \bbC[x, y, z]$ the polynomial ring in three variables.
Let the generators $s, t$ act by permuting the generators $x, y, z$, with $s$ swapping $x$ and $y$, and $t$ swapping $y$ and $z$.
Two examples of morphisms in the category $\Omega(W, Q, \varphi)_{\oplus}$ are
\begin{align*}
    \varphi \colon (st, s, st) \to (st, s), && \psi \colon (t, s) \to (s, st)\\
    \varphi = \begin{blockarray}{cccc}
        & st & s & st \\
        \begin{block}{c(ccc)}
            st & x &  & z \\
            s & & y &  \\
        \end{block}
    \end{blockarray}
    &&
    \psi = \begin{blockarray}{ccc}
        & t & s \\
        \begin{block}{c(cc)}
            s & & x + z  \\
            st & \\
        \end{block}
    \end{blockarray}
\end{align*}
Note that all of the blank entries are zero, and are forced to be zero by definition of the underlying category $\Omega(W, Q, \varphi)$.
The tensor product of these matrices is taken in the usual way (the Kronecker product of matrices), albeit with a twist applied to the matrix on the right depending on what row/column we are in on the left:
\begin{align*}
    \varphi \otimes \psi
    &=
    \begin{blockarray}{cccc}
        & st & s & st \\
        \begin{block}{c(ccc)}
            st & x(st \cdot \psi)  &  & z (st \cdot \psi) \\
            s &  & y(s \cdot \psi) &  \\
        \end{block}
    \end{blockarray}
    \\
    &=
    \begin{blockarray}{ccccccc}
        & s & sts & st & \id & s & sts \\
        \begin{block}{c(cccccc)}
            sts & & x(x + y) & & &  & z(x + y)\\
            ts & & & & & & \\
            \id & & & & y(y + z) & & \\
            t & & & & & & \\
        \end{block}
    \end{blockarray}
\end{align*}
Note finally that $\phi$ and $\psi$ cannot be composed in either direction, since they have distinct domains and codomains.

\subsection{Localisation for the Hecke category}

From this point onwards, we will assume that the realisation $(X, X^\vee)$ is \textit{free} over $\bbZ$, meaning that the simple roots $\Delta \subseteq X$ are linearly independent, so that their $W$-orbit is a root system $\Phi$.
This is not so important for the Hecke category (it does make some arguments simpler), but it is essential for treating localisation in the antispherical case.

Let $(X^\vee, X, \Delta^\vee, \Delta)$ be a free realisation of the Coxeter system $(W, S)$ over the integers, and $\bbk$ a field or complete local ring of arbitrary characteristic.
The Hecke category $\cH_\bbk$ is a linear category whose morphism spaces are $(R, R)$-bimodules where $R$ is the polynomial ring $R = \Sym_\bbk(X)$.
Let $Q$ be $R$ localised at the multiplicative set generated by the roots\footnote{
    We only need to invert the roots, but it would also be fine to take $Q$ to be the field of fractions of $R$ here, and indeed this is what we do in computer algebra software.
} $\Phi$ (in the notation of \cref{section:localisation-antispherical}, we might also write $Q = R_{(\varnothing)}$).
The results of \cite{eliasSoergelCalculus2016} imply that the \defn{localised category} $Q \otimes_{R} \cH_\bbk$ is a semisimple category whose simple objects are in bijection with $W$, with endomorphism ring $Q$ for each simple object.
The standard category $\Std_Q$ defined above (with $\varphi$ the representation of $W$ on $Q$ induced by functorality of $\Sym$ and localisation) is an abstract model for $Q \otimes_{R} \cH_\bbk$, and we will work with the standard category from now on.

\begin{remark}
    When localising a ring, one must be careful to make sure that the set of inverted elements does not contain zero: why does the set of roots $\Phi$ not contain zero after base change?
    By definition, a root is of the form $w \alpha_s$ for some $w \in W$ and $s \in S$, and after base change the action of $w$ is still invertible, and the simple roots $\alpha_s$ are still nonzero.
    Hence $w \alpha_s \neq 0$ inside $X_\bbk$.
    However, previously distinct roots may become equal after base change, and indeed this must occur for cardinality reasons if $\abs{\Phi} \geq p^{\rank X}$, and of course always occurs when $|\Phi| = \infty$.
\end{remark}

There is a monoidal \defn{localisation functor} $\Lambda \colon \cH_\bbk \to \Std_Q$ sending the generating object $B_s$ of $\cH_{BS}$ to the object $(\id, s)$, and the generating morphisms of $\cH_{BS}$ (up-dots, down-dots, trivalent vertices, and $2m_{st}$-valent vertices) to matrices over $Q$: these matrices are calculated explicitly in \cite{eliasLocalizedCalculusHecke2020}.
This gives a right action of $\cH_{BS}$ on $\Std_Q$, where $X \cdot Y = X \otimes \Lambda(Y)$ for $X \in \Std_Q$ and $Y \in \cH_{BS}$.

The main idea of the algorithm is then to follow the light leaves construction in order to calculate, for each element $x \in W$:
\begin{enumerate}
    \item A \defn{standard object} $T_x \in \Std_Q$ such that $\Lambda(B_x) = T_x$,
    \item For all $y \leq x$, the sets $\TLLup(T_y, T_x) := \Lambda(\TLLup(B_y, B_x))$ and $\TLLdown(T_x, T_y) := \Lambda(\TLLdown(B_x, B_y))$ of truncated light leaves inside of $\Hom_{\Std_Q}(T_y, T_x)$ and $\Hom_{\Std_Q}(T_x, T_y)$ respectively\footnote{
        Actually, we need only keep a subset of these elements which descend to a basis in the local quotient: see \cref{remark:efficiency}.
}.
\end{enumerate}
For example we start with $T_{\id} = \mathbf{1} \in \Std_Q$, the monoidal unit, and up and down light leaves simply the identity.
For each $s \in S$, the object $T_s = \Lambda(\mathbf{1}) \cdot B_s$ (we know this object will always be indecomposable, no matter the characteristic).
Then $\TLLup(T_{\id}, T_s) = \TLLup(T_{\id}, T_{\id}) \cdot \bsunit$ and $\TLLdown(T_s, T_{\id}) = \TLLdown(T_{\id}, T_{\id}) \cdot \bscounit$, according to the light leaves construction.

Things get more complicated later in the process, where we cannot assume that ${^\bbk B_x} \cdot B_s$ is indecomposable in $\cH_{\bbk}$.
In this case we instead have generated sets $\TLLup(T_y, T_x \cdot B_s)$ and $\TLLdown(T_x \cdot B_s, T_y)$ which capture the corresponding sets $\TLLup({^\bbk B_y}, {^\bbk B_x} \cdot B_s)$ and $\TLLdown({^\bbk B_x} \cdot B_s, {^\bbk B_y})$.
Using the theory of intersection forms, these sets of morphisms can be used to find a direct summand $T_{xs} \subseteq T_x \cdot B_s$ which represents the indecomposable ${^\bbk B}_{xs}$, then the calculated morphisms can be trimmed (along the projection and inclusion defined by the direct summand) into $\TLLup(T_y, T_{xs})$ and $\TLLdown(T_{xs}, T_y)$.

\begin{remark}
    The standard object $T_x$ is an ordered sequence of elements of $W$, and the multiset of those elements is invariant under isomorphism in $\Std_Q$.
    If $p = 0$ then this multiset gives the multiplicities $h_{y, x}(1)$ of Kazhdan-Lusztig polynomials at $v = 1$, and analagously for $p > 0$ the multiset gives the multiplicities ${^p h_{y, x}}(1)$ of the $p$-Kazhdan-Lusztig polynomials at $v = 1$.
\end{remark}

This brief description is just a guide: the actual algorithm is stated in \cref{section:algorithm}.

\subsection{Localisation for the antispherical category} \label{section:localisation-antispherical}

Let $(X^\vee, X, \Delta^\vee, \Delta)$ be a realisation of the Coxeter system $(W, S)$ over the integers, $\bbk$ a field or complete local ring of arbitrary characteristic, and $I \subseteq S$ a subset of the simple reflections.
Localisation for the antispherical category is necessarily different than from the Hecke category: no matter how we define our localisation functor $\Lambda_I$, we should for example have $\Lambda_I(\mathbf{1}) \cdot B_s = 0$ for $s \in I$.
We need the parabolic roots $\Phi_I$ to act by zero, while we still need many other roots inverted to guarantee semisimplicity of the resulting category, so we want to invert the set of non-parabolic roots $\Phi \setminus \Phi_I$ in $R$, then take the quotient by $\innprod{\alpha_s \mid s \in I}$.

In characteristic zero, a root $\beta = \sum_{s \in S} b_s \alpha_s$ is in $\Phi_I$ if and only if $b_s = 0$ for all $s \notin I$, therefore every $\beta \in \Phi \setminus \Phi_I$ remains nonzero after taking the quotient by $\innprod{\alpha_s \mid s \in S}$.
After reduction modulo $p$ however, a non-parabolic root $\beta \in \Phi_I \setminus \Phi$ may ``look'' parabolic in the sense that it becomes zero in the quotient.
Whenever this happens, we have inverted an element $\beta$ and then set it to zero, and so the ring $Q_I$ just defined is zero: and our process of localisation will not work.
The remedy to this problem is to instead work with a $p$-modular system, letting the realisation live over a local ring of characteristic $0$ whose residue field has characteristic $p$.

Let $\bbO$ be a local ring of characteristic zero, $\bbK$ its field of fractions (also of characteristic zero), and $k$ its residue field of characteristic $p > 0$: then $(\bbO, \bbK, k)$ is called a \defn{$p$-modular system}.
We will most often use $(\bbZ_{(p)}, \bbQ, \bbF_p)$\footnote{
    This choice is because it allows simple exact computation in MAGMA.
}, where $\bbZ_{(p)} \subseteq \bbQ$ is the subring of the rationals consisting of fractions whose denominator is not divisible by $p$.
Base change the realisation $(X, X^\vee, \Delta, \Delta^\vee)$ to the local ring $\bbO$, and define the Hecke category $\cH_\bbO$ using the local ring $\bbO$, so we are working over the polynomial ring $R = \Sym_\bbO(X)$.
The indecomposable objects ${^\bbO}B_w$ of $\cH_\bbO$ base change over $k$ to the indecomposable objects ${^k} B_w$, hence we can use the former to make arguments about the latter.

Since $\bbO$ has characteristic zero, we can define the localised ring $R_{(I)} = \bbK \otimes R[\beta^{-1} \mid \beta \in \Phi \setminus \Phi_I]$, and the resulting quotient $Q_I = R_{(I)} / \innprod{\alpha_s \mid s \in I}$ is nonzero.
We modify the definition of the groupoid category a little: the endomorphism space $\End(r_x)$ for $x \in {^I W}$ will be a $(Q_I, R)$-bimodule, with the map ${^I \varphi} \colon {^I W} \to \Hom(R, Q_I)$ being the action of $x$ on $R$, followed by the natural map $R \to Q_I$.
Let ${^I \Std_{Q_I}}$ be the groupoid category $\Omega({^I W}, R, Q_I, {^I \varphi})$.

The antispherical category ${^I \Std_{Q_I}}$ admits a right action of $\cH_{BS}$, and an \defn{antispherical localisation functor} $\Lambda_I \colon {^I \cN_\bbO} \to {^I \Std_{Q_I}}$ commuting with this action.
The right action of $B_s$ on a simple object $x$ is
\begin{equation} \label{equation:bs-antispherical-action}
    (x) \cdot B_s = \begin{cases}
        (x, xs) & \text{if } xs \in {^I W}, \\
        0 & \text{otherwise}.
    \end{cases}
\end{equation}
The action of a generating morphism $\psi \colon \underline{w} \to \underline{x} \in \cH_{BS}$ can be deduced from the data of $\underline{w}$, $\underline{x}$, and the image $\Lambda(\psi) \in \Std_Q$ of the generating morphism in the full standard category, which we show in the example below.
(Note: we really do need to remember $\underline{w}$ and $\underline{x}$: just knowing the domain and codomain of $\Lambda(\psi)$ is not enough).

\textit{Example:} Let $W = \innprod{s, t \mid s^2 = t^2 = (st)^3 = 1}$ be the symmetric group on 3 letters, and $\alpha, \beta \in X$ the simple roots corresponding to $s$ and $t$ respectively.
The image of the $2m_{st}$-valent vertex $\psi \colon (s, t, s) \to (t, s, t)$ is
\begin{equation}
    \Lambda(\psi)
    =
    \begin{blockarray}{ccccccccc}
    & \id                          & s & t                            & ts & s & \id                          & st & sts\\
        \begin{block}{c(cccccccc)}
\id & \frac{\alpha + \beta}{\beta} &   &                              &    &   & \frac{\alpha + \beta}{\beta} &    & \\
t   &                              &   & \frac{\alpha + \beta}{\beta} &    &   &                              &    & \\
s   &                              & 1 &                              &    & 1 &                              &    & \\
st  &                              &   &                              &    &   &                              & 1  & \\
t   &                              &   & -\frac{\alpha}{\beta}        &    &   &                              &    & \\
\id & -\frac{\alpha}{\beta}        &   &                              &    &   & -\frac{\alpha}{\beta}        &    & \\
ts  &                              &   &                              & 1  &   &                              &    & \\
tst &                              &   &                              &    &   &                              &    & 1 \\
        \end{block}
    \end{blockarray}
\end{equation}
Now suppose that we take the parabolic subgroup defined by $I = \{t\}$, so that ${^I W} = \{\id, s, st\}$.
We will compute the action of $\psi \in \cH_{BS}$ on the identity map $1_s \in \Std_I$.
Firstly we act on the domain and codomain to get (recall \cref{equation:bs-antispherical-action}: we are crossing out terms which are zero)
\begin{align}
    (s) &\cdot B_s \cdot B_t \cdot B_s                                                       & (s) &\cdot B_t \cdot B_s \cdot B_t  \\
                         & = (s, \id) \cdot B_t \cdot B_s                               &                      & = (s, st) \cdot B_s \cdot B_t \\
                         & = (s, st, \cancel{\id, t}) \cdot B_s                &                      & = (s, \id, \cancel{st, sts}) \cdot B_t \\
                         & = (s, \id, \cancel{st, sts}, \cancel{\id, s, t, ts})  &                      & = (s, st, \cancel{\id, t}, \cancel{st, s, sts, ts}),
\end{align}
so $1_s \cdot \psi$ is a map from $(s, \id)$ to $(s, st)$ in $\cat{Std}_I$.
From our calculations of the action on objects, only the first two entries out of 8 ``survive'' in both the domain and codomain, so the matrix of $1_s \cdot \psi$ is formed by taking the top-left $2 \times 2$ submatrix of $\Lambda(\psi)$, twisting it by $s$, and then performing the specialisation map $\pi_I \colon R_{(I)} \surjto Q_I$.
Note that $s(\alpha) = -\alpha$ and $s(\beta) = \alpha + \beta$, so
\begin{equation}
    s \cdot \Lambda(\psi)_{2 \times 2}
    =
    \begin{blockarray}{ccc}
               & s & \id \\
        \begin{block}{c(cc)}
        s                                                                                              & s \cdot \left(\frac{\alpha + \beta}{\beta}\right)                    & \\
        st                                                                                             &                                                                      & \\
        \end{block}
    \end{blockarray}
    =
    \begin{pmatrix} \frac{\beta}{\alpha + \beta} & \\ & & \end{pmatrix}
    \xsurjto{\pi_I}
    \begin{blockarray}{ccc}
               & s & \id \\
        \begin{block}{c(cc)}
        s                                                                                              & \frac{0}{\alpha}                                                     & \\
        st                                                                                             &                                                                      & \\
        \end{block}
    \end{blockarray}
\end{equation}
Hence $1_s \cdot \psi \colon (s,\id) \to (s, st)$ is the zero map.
Notice that something nice has happened above: even though $\beta = 0$ in $Q_I$ because $\beta \in \Phi_I$, something has conspired to ensure that $\beta$ does not appear in the denominators of our resulting twisted submatrix $s \cdot \Lambda(\psi)_{2 \times 2}$, so that we can indeed set $\beta \mapsto 0$.
For the procedure outlined here to make sense, we need to show that every element in the twisted submatrix lands in the subring $R_{(I)}$ of $Q$.

\begin{itemize}
    \item
        If $\psi$ is a one-colour generator (the unit, counit, multiplication, or comultiplication), then its domain and codomain are products of $B_s$ for some $s$, and the only root appearing in the denominators of $\Lambda(\psi)$ is $\alpha_s$.
        Let $x \in {^I W}$, then if $xs \in {^I W}$ we have by the parabolic property that $x(\alpha_s) \notin \Phi_I$ so all entries of $x \cdot \Lambda(\psi)$ make sense, and $1_x \cdot \Psi$ is defined.
        If $xs \notin {^I W}$ then $x \cdot B_s = 0$ so there is nothing to check.
    \item
        The only two-colour generator that is needed is the ``spider'' morphism (the $2 m_{s,t}$-valent vertex, with one end turned around, so it is a morphism from $B_s B_t B_s \cdots \to \mathbf{1}$).
        The matrix $\Lambda(\psi)$ in this case is a polynomial: the product of all positive roots in the dihedral subgroup generated by $\alpha_s$ and $\alpha_t$.
        So there is no zero division.
\end{itemize}

\begin{remark}
    The reader will notice that this argument was rather subtle and used explicit knowledge of the localisations of the generating morphisms of $\cH_{BS}$.
    It seems somewhat miraculous that this works, and we are currently lacking a more conceptual explanation.
\end{remark}

\subsection{The local quotient and intersection form} \label{section:local-intersection-form}

We have mentioned the local quotient functor before in the context of the Hecke category: given some $w \in W$, the \defn{local quotient} $\cH_{\bbk}^{\not < w} = \cH_\bbk / \cH_{\bbk, < w}$ is the quotient of the Hecke category by morphisms factoring through the additive category generated by the ${^\bbk B}_x$ for $x < w$, and the \defn{local quotient functor} is the natural functor $\cH_\bbk \surjto \cH_{\bbk}^{\not < w}$.
Note that this functor is not monoidal, as the quotient category is not monoidal.
In the local quotient category $\cH_{\bbk}^{\not < w}$ the endomorphism ring of ${^\bbk B}_w$ becomes isomorphic to $R_\bbk$, which is plain to see from the light leaves construction: the sets $\TLLdown({^\bbk B}_w, {^\bbk B}_w)$ and $\TLLup({^\bbk B}_w, {^\bbk B}_w)$ are each free of rank 1, spanned by the identity morphism, while all lower leaves become zero in the quotient.

In both algorithms we present, we use the local quotient functor to count multiplicities of indecomposables inside other objects.
Let $X$ be any object of the Hecke category $\cH_{\bbk}$ which we expect might contain one or more copies of the indecomposable ${^\bbk B_w}$.
The \defn{local intersection form} of $X$ at $w \in W$ is the $R$-valued bilinear form
\begin{align}
    I_{X, w} \colon \Hom_{\not < w}(X, {^\bbk B_w}) \times \Hom_{\not < w}({^\bbk B_w}, X) &\to \End_{\not < w}({^\bbk B_w}) \cong R_\bbk, \\
    \left({^\bbk B_w} \xleftarrow{f} X\right) \times \left(X \xleftarrow{g} {^\bbk B_w} \right) &\mapsto f \circ g.
\end{align}
By Corollary 11.75 of \cite{eliasIntroductionSoergelBimodules2020}, the graded rank of $I_{X, w}$ is the graded multiplicity of the indecomposable object ${^\bbk B_w}$ inside the larger object $X$.

\begin{remark} \label{remark:efficiency}
    In the above equation, we need only consider a set of those $f \in \TLLdown(X, {^\bbk B_w})$ which descend to a basis in the quotient space $\Hom_{\not <w}(X, {^\bbk B_w})$, and likewise for the $g$.
    As it turns out, throughout the algorithm we only ever need to keep track of a subset of $\TLLdown$ and $\TLLup$ which descend to a basis in the appropriate quotient.
    This is crucial to the efficiency of the algorithm.
\end{remark}

There is an analogue of the quotient category $\cH_{\bbk}^{\not <w}$ in the localised case.
In the notation of \cref{section:standard-category}, let $\Std_Q^{\not < w}$ be the category $\Omega(Q, W^{\not <w}, \varphi)_{\oplus}$, then there is a natural quotient functor $\cat{Loc}^{\not <w} \colon \Std_Q \surjto \Std_Q^{\not < w}$ which sends the simple object $r_x$ to zero if $x < w$.
(Again, this functor is not monoidal because its target is not monoidal).
In the concrete interpretation of a morphism being a framed matrix with entries in $Q$, the action of $\cat{Loc}^{\not <w}$ is to extract the submatrix whose framed entries are all $\not <w$.
For a standard object $T_w$ which represents ${^\bbk B}_w$, every $r_x$ summand has $x \leq w$, and there is a unique summand (using our convention for the tensor product order, the last element of the sequence) which is $r_w$.
So the local quotient functor $\cat{Loc}^{\not <w}$ is particularly simple on morphisms into or out of $T_w$: extract the last row or column of the matrix respectively, deleting any other entries whose framing elements are necessarily $< w$.


\section{The algorithm} \label{section:algorithm}

Throughout this section, let $(\bbO, \bbK, k)$ be a $p$-modular system, our algorithm uses $(\bbZ_{(p)}, \bbQ, \bbF_p)$.
Our goal is to build a model for the Hecke (or antispherical) category inside $\Std_Q$ (or ${^I \Std}_{Q_I}$).
Since most things will be the same, we will just state the algorithm for the Hecke category, and point out any places that the story changes for the antispherical category --- the main change being the replacement of the standard category and its modified right action by $\cH_{BS}$.

The Hecke category $\cH_\bbO$ is an additive Krull-Schmidt\footnote{
    A \defn{Krull-Schmidt category} is an additive category where every object decomposes into a finite direct sum of indecomposable objects, and indecomposable objects are precisely those objects with local endomorphism rings.
    In such a category the \defn{Krull-Schmidt property} holds: when an object is written as a direct sum of indecomposables, the multiplicities of the (isomorphism classes of) indecomposables does not depend on the decomposition chosen.
    In other words, each object has finite length, and the multiplicities of indecomposable objects inside arbitrary objects is well-defined.
    (Not every additive category having the Krull-Schmidt property is a Krull-Schmidt category: for example $\bbZ$ is indecomposable in the category of finitely generated abelian groups, but $\End_\bbZ(\bbZ) \cong \bbZ$ is not a local ring).
} category, linear over a ring $R = \Sym_\bbO^\bullet(X)$, with indecomposable objects $\set{{^\bbO B}_w \mid w \in W}$ parameterised by $w \in W$ up to isomorphism and grading shift.
Hence knowledge of the category $\cH_\bbO$ is the data of the sets $\Hom_{\cH_\bbO}({^\bbO B_x}, {^\bbO B_y})$ for all $x, y \in W$, together with the $R$-module structure on the hom-set, and a rule for composing compatible morphisms.
We will provide such a description of the Hecke category, with two major simplifications:
\begin{enumerate}
    \item
        Since we have a faithful functor to $\Std_Q$, where composition of morphisms is simply matrix multiplication, we only need to specify the image of the hom spaces and composition will come for free.
        Better yet, the category $\Std_Q$ carries its own left $R$-action via the natural map $R \to Q$, so the $R$-action comes for free; furthermore in order to specify the image of a hom space, we need only remember an $R$-spanning set.
    \item
        Because of the double leaves construction, rather than remembering full hom spaces we may instead just keep track of the smaller spaces $\TLLdown$ and $\TLLup$, as the hom-spaces may be reconstructed from them.
\end{enumerate}

We now introduce a little language; by a \textit{model} for the category $\cH_\bbO$, we mean the following data:
\begin{enumerate}
    \item A set $\set{T_x \mid x \in W}$ of objects in $\Std_Q$,
    \item For each $x \leq y$, finite sets $\TLLup(T_x, T_y) \subseteq \Hom_{\Std_Q}(T_x, T_y)$ and $\TLLdown(T_y, T_x) \subseteq \Hom_{\Std_Q}(T_y, T_x)$.
    \item For each element $w \in W$ and right descent $s \in \cR(w)$, an \textit{inclusion map} $i_{w, s} \colon T_w \to T_{ws} \cdot B_s$ and a \textit{projection map} $p_{w, s} \colon T_{ws} \cdot B_s \to T_w$.
\end{enumerate}
The finite set $\TLLup(T_x, T_y)$ is considered to be a set of vectors spanning an $R$-lattice.
For example, redundant (linearly dependent) vectors may be removed, provided the linear dependencies are over $R$ (not over $Q$).
The inclusions and projections $i_{w, s}$ and $p_{w, s}$ are not needed for the data of the category, but they will be needed for using the right action of $\cH_{BS}$ in order to iterate the (truncated) light leaves construction.
Morally they are the inclusion and projection maps for the unique indecomposable summand $B_w$ appearing in $B_{ws} \cdot B_s$.
Keeping track of them allows us to work directly with indecomposable objects rather than Bott-Samelson objects.

\subsection{Overview of the algorithm}

Here we provide a broad overview of how the different parts of the algorithm fit together, and fix notation going forward.
We have tried to give suggestive names to each step of the algorithm, and expand upon these steps in subsections.
The algorithm processes group elements in an order compatible with the Bruhat order; from this point forward fix a $w \in {^I W}$, and suppose that we have already processed all elements less than $w$ in the Bruhat order.

For $w \in {^I W}$ and each right descent $s \in \cR(w)$, we produce some preliminary data:
\begin{enumerate}[label=(\alph*)]
    \item \textit{Translate morphisms:} The sets $\TLLup(-, T_{ws} \cdot B_s)$ and $\TLLdown(T_{ws} \cdot B_s, -)$ modelling the truncated light leaves in and out of $B_{ws} B_s$ are computed from data of $\TLLup(-, T_{ws})$, $\TLLdown(T_{ws}, -)$, the inclusions $i_{x, s}$ and projections $p_{x, s}$ for all smaller elements $x < w$, and the action of $\cH_{BS}$ on $\Std_I$.
    This step roughly follows the construction of light leaves in the Bott-Samelson category.
    \item \textit{Reduce idempotent:} The truncated light leaves in and out of $T_{ws} \cdot B_s$ are used to find any summands of $T_{ws} \cdot B_s$ isomorphic to $T_x$ for $x < w$, and their corresponding idempotents. Subtracting all of these from the identity map on $T_{ws} \cdot B_s$ yields an idempotent $E_w^s \in \End_{\Std}(T_{ws} \cdot B_s)$ modelling the idempotent in $\End_{\cH}(B_{ws} B_s)$ projecting to the top summand $B_w$.
    \item \textit{Split idempotent:} The idempotent $E_w^s$ is split into an isomorphism onto its image, followed by an inclusion back into $T_{ws} \cdot B_s$. More precisely, we find a subobject $T_w^s \subseteq T_{ws} \cdot B_s$ which is $\cat{Std}_I$-isomorphic to $\im E_w^s$, and a splitting
        \begin{align}
            i_w^s &\colon T_w^s \to T_{ws} \cdot B_s &
            p_w^s &\colon T_{ws} \cdot B_s \to T_w^s
        \end{align}
        such that $p_w^s \circ i_w^s = \id_{T_w^s}$ and $i_w^s \circ p_w^s = E_w^s$.
\end{enumerate}
There may be many right descents of $w$: an arbitrary one $s_0 \in \cR(w)$ is chosen to define the \defn{standard sequence} $T_w := T_w^{s_0}$.
There are then two more steps to take before we are finished:
\begin{enumerate}
    \item \textit{Braiding:} For each $s \in \cR(w)$, the \defn{inclusion} and \defn{projection} maps
        \begin{align}
            i_{w, s} &\colon T_w \to T_{ws} \cdot B_s &
            p_{w, s} &\colon T_{ws} \cdot B_s \to T_w &
            \text{(defined when } ws < w \text{)}
        \end{align}
        are defined by either postcomposing or precomposing with an isomorphism between $T_x^s$ and $T_x = T_x^{s_0}$.
        For $s = s_0$ we have $i_{w, s_0} = i_w^{s_0}$ and $p_{w, s_0} = p_w^{s_0}$; otherwise the $2m_{s, s_0}$-valent vertex is used is produce the isomorphism.
    \item \textit{Morphism truncation:} For $s = s_0$, the morphisms $\TLLup(T_{ws} \cdot B_s)$ are postcomposed (``truncated'') by $p_{w, s}$ to yield a set of maps into $T_w$, and similarly for $\TLLdown(T_{ws} \cdot B_s)$ with $i_{w, s}$.
        These ``truncated'' maps may have become linearly dependent over $R$.
    \item
        \textit{Morphism pruning:} A subset of linearly independent maps is selected from the truncated maps to serve as the new set of truncated light leaves $\TLLup(T_w)$ and $\TLLdown(T_w)$.
        As mentioned in \cref{remark:efficiency}, for each $x \leq w$ we only actually keep a subset of $\TLLup(T_x, T_w)$ corresponding to a basis in the local quotient $\cH_{\bbO}^{\not < x}$.
\end{enumerate}
The data $(w, T_w, (i_{w, s})_{s \in \cR(w)}, (p_{w, s})_{s \in \cR(w)}, \TLLup(T_w), \TLLdown(T_w))$ is the output of the algorithm at each step, and used in calculations for the subsequent larger elements. A few notes:
\begin{itemize}
    \item The inductive process starts with the identity element $\id \in W$, for which the data is straightforward: we have the standard sequence $T_{\id} = (\id)$, no projections or inclusions since there are no right descents, and the only truncated light leaf is the identity map so
    \[\TLLup(T_{\id}, T_{\id}) = \TLLdown(T_{\id}, T_{\id}) = \{\id_{T_{\id}}\}.\]
    \item The action of the one-colour generators of $\cH_{BS}$ is used during the \textit{translate morphisms} step, and the action of the $2m_{st}$-valent vertex is used during the \textit{braiding} step.
        These are the only places where the action of the category is used. In particular, the $2m_{st}$-valent vertex has been encoded into the set of inclusions and projections, which means that it can be avoided entirely during the \textit{translate morphisms} step, contrary to the usual construction of light leaves appearing in \cite{eliasIntroductionSoergelBimodules2020}.
    \item The only significant differences between the characteristic zero and characteristic $p > 0$ cases arise in the \textit{reduce idempotent} and \textit{morphism pruning} steps, where we need to be careful to work with an $R$-lattice inside a $Q$-vector space, $R$ being a graded polynomial ring over $\bbO$ and $Q$ being a rational function ring over $\bbK$.
        The rest of the algorithm is oblivious to the chosen characteristic.
    \item The \textit{reduce idempotent}, \textit{split idempotent}, and \textit{morphism pruning} steps are highly noncanonical, and in general there will be many choices that yield isomorphic models of categories.
    Together with the fact that the truncated light leaves themselves are highly noncanonical, this may be a source of trouble in the long run (see \cref{section:limitations}).
    \item We are very interested in the self-dual coefficients ${^p m_{x, w}} \in \bbZ_{\geq 0}[v^{\pm 1}]$ expressing the $p$-canonical basis element ${^p b_w}$ in terms of the canonical basis elements ${^p b_w} = \sum_{x \leq w} {^p m_{x, w}} b_w$.
        These coefficients can be recovered as the graded rank of the set $\TLLup(T_x, T_w)$, but also are produced naturally during the computation of (b).
\end{itemize}

\subsection{Translate morphisms}

Let $w \in {^I W}$ with a right descent $s$ so that $w > ws$, and let $x \leq ws$. If $x < xs$ we define two homomorphisms
\begin{align}
    U_0(-, s) &\colon \Hom_{\cat{Std}}(T_{ws}, T_x) \to \Hom_{\cat{Std}}(T_{ws} \cdot B_s, T_x), \\
    U_1(-, s) &\colon \Hom_{\cat{Std}}(T_{ws}, T_x) \to \Hom_{\cat{Std}}(T_{ws} \cdot B_s, T_{xs}),
\end{align}
which mimic the light leaves construction. For $\alpha \colon T_{ws} \to T_x$ a previously constructed truncated light leaf down, the morphisms $U_0(\alpha, s)$ and $U_1(\alpha, s)$ are the top and bottom paths respectively in the following diagram:
\begin{equation}\label{cd:lloutU}
    \begin{tikzcd}[row sep=0em]
                                                       &                                                                              & T_x\\
        T_{ws} \cdot B_s \arrow[r, "\alpha \cdot B_s"] & T_x \cdot B_s \arrow[ru, "T_x \cdot \bscounit"] \arrow[rd, "{p_{xs, s}}"] \\
                                                       &                                                                              & T_{xs}
    \end{tikzcd}
\end{equation}
In the special case $x = ws$, we have $xs = w$ and so $U_1(\alpha, s)$ does not make sense since we have not yet defined $p_{w, s}$. In this case, set $U_1(\alpha, s) = \alpha \cdot B_s$.

If instead we have $x > xs$, define two homomorphisms
\begin{align}
    D_0(-, s) &\colon \Hom_{\cat{Std}}(T_{ws}, T_x) \to \Hom_{\cat{Std}}(T_{ws} \cdot B_s, T_x), \\
    D_1(-, s) &\colon \Hom_{\cat{Std}}(T_{ws}, T_x) \to \Hom_{\cat{Std}}(T_{ws} \cdot B_s, T_{xs}),
\end{align}
where $D_0(\alpha, s)$ and $D_1(\alpha, s)$ are the top and bottom paths respectively in the following diagram:
\begin{equation}\label{cd:lloutD}
    \begin{tikzcd}[row sep=0em]
        & & & & T_x \\
        T_{ws} \cdot B_s \arrow[r, "\alpha \cdot B_s"] &
        T_x \cdot B_s \arrow[r, "{i_{x, s}}"] &
        T_{xs} \cdot B_s \cdot B_s \arrow[r, "T_{xs} \cdot \bsmult"] &
        T_{xs} \cdot B_s \arrow[ur, "{p_{x, s}}"] \arrow[dr, "T_{xs} \cdot \bscounit"]
        \\
        & & & & T_{xs}
    \end{tikzcd}
\end{equation}

Taking all truncated light leaves $\TLLdown(T_{ws}, -)$ and applying $U_0, U_1, D_0, D_1$ (whichever two are valid for the given morphism) generates the new set $\TLLdown(T_{ws} \cdot B_s, -)$. Explicitly, if we relabel $x$ and $xs$ so that $x > xs$ we have
\begin{align}
    \TLLdown(T_{ws} \cdot B_s, T_x) = \bigsqcup_{\beta \in \TLLdown(T_{ws}, T_{xs})} U_1(\beta, s) \sqcup \bigsqcup_{\alpha \in \TLLdown(T_{ws}, T_x)} D_0(\alpha, s) && (x > xs)\\
    \TLLdown(T_{ws} \cdot B_s, T_{xs}) = \bigsqcup_{\beta \in \TLLdown(T_{ws}, T_{xs})} U_0(\beta, s) \sqcup \bigsqcup_{\alpha \in \TLLdown(T_{ws}, T_x)} D_1(\alpha, s) && (x > xs) \\
\end{align}
There is one last case that may happen in the antispherical category, which is that $x \leq ws$ but $xs \notin {^I W}$.
In this case we simply set $\TLLdown(T_w \cdot B_s, T_{xs}) = \varnothing$.

The construction of $\TLLup(T_{ws} \cdot B_s)$ from $\TLLup(T_{ws})$ is analagous.
The diagrams for $U_0, U_1, D_0$, and $D_1$ are read backwards, replacing any projections with their corresponding inclusions, and actions of generators in $\cH_{BS}$ by their upside-down versions.

\subsection{Reduce idempotent}

During the \textit{reduce idempotents} step, the only data we need are the $\TLLup(T_{ws} \cdot B_w)$ and $\TLLdown(T_{ws} \cdot B_w)$ morphisms.
To begin with, initialise the idempotent $E$ to the identity morphism on $T_{ws} \cdot B_w$.

For each element $x < w$ and degree $d \in \bbZ$, define the \textit{local intersection form} to be the matrix
\[M = [\overline{\text{out}_i \circ E \circ \text{in}_j}]_{i,j},\]
where $\text{in}_j$ ranges over $\TLLup^d(T_x, T_{ws} \cdot B_s)$ and $\text{out}_i$ ranges over $\TLLdown^{-d}(T_{ws} \cdot B_w, T_x)$, and the overline $\overline{f}$ denotes the projection of $f$ to the category $\Std^{\not <x}$.
In particular, as the matrix $\overline{\text{in}_j}$ has one column and the matrix $\overline{\text{out}_i}$ has one row, the intersection form $M$ is a matrix with entries in $Q$.
In fact it has entries in $\bbO$, since we are modelling an $R_\bbO$-linear category.

The rank of the matrix $M$ modulo $p$ (i.e. the rank of $M \otimes_{\bbO} k$) gives the number of times the indecomposable $B_x$ appears in the idempotent $E$.
If $M$ has nonzero rank, then let $(i, j)$ be any pair such that $\overline{\text{out}_i \circ E \circ \text{in}_j}$ is invertible in $\bbO$, it then follows from an idempotent lifting argument that $\text{out}_i \circ E \circ \text{in}_j$ is invertible, so that
\begin{equation}
    E \circ \text{in}_j \circ (\text{in}_j \circ E \circ \text{out}_i)^{-1} \circ \text{out}_j \circ E
\end{equation}
is an idempotent summand in $E$.
This summand is subtracted off $E$, and the intersection form is re-calculated, repeating the above process until it has zero rank.
When we finish, $E$ is an idempotent morphism in $\End_{\Std_Q}(T_{ws} \cdot B_s)$ projecting to the top term, the image of $B_w$.

\begin{remark}
    Since the $p$-canonical basis element ${^p b_{ws}}$ is already known, the character of $B_{ws} \cdot B_s$ can be calculated in the Hecke algebra as ${^p b_{ws}}b_s$.
    We can then calculate the $p$-canonical basis element ${^p b_w}$ by subtracting off the sum of the ${^p b_x}$ for each $x < w$ removed during the reduce idempotents step.
    This gives another method other than graded rank to arrive at ${^p b_w}$, and can serve as a useful integrity check throughout the algorithm.
\end{remark}

\subsection{Split idempotent}

Any idempotent matrix may be decomposed into the composition of projection onto its image, followed by inclusion of its image back into the ambient space.
In general there are many ways of doing this.
Our method of splitting the idempotent $E_{w}^s$ does not rely on any special structure of the category; it is essentially just decomposing an idempotent matrix of rank $r$ into a product of two rank $r$ matrices, while keeping track of the framing data which associates rows and columns with elements of $W$.
There are many ways that this could be done, we write one down for concreteness.

Let $E \in \End_{\Std}((x_1, \ldots, x_n))$ be an idempotent, so that $E^2 = E$. Let \[i \colon (x_{\sigma(1)}, \ldots, x_{\sigma(r)}) \to (x_1, \ldots, x_n)\] be the reduced column-echelon form of the matrix $E$ with zero columns discarded, where $\sigma$ is some injective map of $\sigma \colon \{1, \ldots, r\} \to \{1, \ldots, n\}$ and $r$ is the rank of $E$. Let $p \colon (x_1, \ldots, x_n) \to (x_{\sigma(1)}, \ldots, x_{\sigma(r)})$ be any matrix solving the equation $E = ip$.
The fact that $E^2 = E$ implies that $ipip = ip$; since $i$ is injective and $p$ is surjective we may left and right cancel to get $pi = \id$, so indeed $p$ and $i$ split the idempotent $E$.

Applying this process to $E_w^s \in \End_{\Std_Q}(T_{ws} \cdot B_s)$ yields $i_w^s$ and $p_w^s$: the domain of $i_w^s$ (or equivalently the codomain of $p_w^s$) determine the standard sequence $T_w^s$ which represents the object ${^\bbO B}_w$.

\subsection{Braiding}

By now, we have run all of the above steps for all choices of right descent $s \in \cR(w)$.
Suppose that $w$ has two distinct right descents $s, t \in \cR(w)$.
We have generated two different objects $T_w^s$ and $T_w^t$ which model $B_w$.
They are isomorphic but not in general equal, and our goal is to find an isomorphism $\varphi_{s, t} \colon T_w^s \to T_w^t$ between them which comes from an isomorphism in $\cH_\bbO$.
The distinguished isomorphism is based on using the inclusion maps to move from $T_w^s$ to $T_x \cdot (\cdots B_s B_t B_s)$, where $x$ is the shortest element in the right $\{s, t\}$ coset containing $w$, applying the $2 m_{st}$-valent vertex to go to $T_x \cdot (\cdots B_t B_s B_t)$, and then using the projection maps to move back to $T_w^t$.
Based on this description the reader should be able to figure out how to write down the morphism $\varphi_{s, t}$, but we do the case $m_{st} = 3$ for clarity (note that $wsts = x = wtst$)
\begin{equation}
    \begin{tikzcd}[column sep=4em]
        T_w^s \arrow[r, "i_w^s"] \arrow[d, dotted, "\varphi_{s, t}"] & T_{ws} \cdot B_s \arrow[r, "{i_{ws, t} \cdot B_s}"] & T_{wst} \cdot B_t B_s \arrow[r, "{i_{wst, s} \cdot B_t B_s}"] & T_{wsts} \cdot B_s B_t B_s \arrow[d, "{T_x \cdot \text{braid}}"]\\
        T_w^t                                                & T_{wt} \cdot B_t \arrow[l, "p_w^s"]                 & T_{wts} \cdot B_s B_t \arrow[l, "p_{wt,s} \cdot B_t"]         & T_{wtst} \cdot B_t B_s B_t \arrow[l, "p_{wts,t} \cdot B_s B_t"]
    \end{tikzcd}
\end{equation}
The maps $\varphi_{s, t}$ and $\varphi_{t, s}$ are inverse to each other, and in what follows we set $\varphi_{s, s} = \id$.

With the maps $\varphi_{s, t}$ defined, the actual job of the \textit{braiding} step is to define the inclusion and projection maps $i_{w, s}$ and $p_{w, s}$ relative to the chosen descent $s_0$.
They are defined by the commutative diagrams
\begin{equation}
    \begin{tikzcd}
        T_w^s \arrow[r, "i_w^s"] & T_{ws} \cdot B_s \\
        T_w = T_w^{s_0} \arrow[u, "\varphi_{s_0, s}"] \arrow[ur, dotted, "{i_{w, s}}"]
    \end{tikzcd}
    \quad
    \text{and}
    \quad
    \begin{tikzcd}
        T_w^s \arrow[d, "\varphi_{s,s_0}"] & T_{ws} \cdot B_s \arrow[l, "p_w^s"] \arrow[ld, dotted, "{p_{w, s}}"] \\
        T_w = T_w^{s_0}
    \end{tikzcd}
\end{equation}

Once the maps $i_{w, s}$ and $p_{w, s}$ are generated for all $s \in \cR(w)$, the auxiliary maps $i_w^s$ and $p_w^s$ may be forgotten.

\subsection{Morphism truncation}

The purpose of this step is to use the projection $p_{w, s_0}$ and inclusion $i_{w, s_0}$ to narrow the set of light leaves $\TLLup(T_{ws_0} \cdot B_{s_0})$ and $\TLLdown(T_{ws_0} \cdot B_{s_0})$ so that they give an $R$-span for the light leaves in and out of the modelled object $T_w = T_w^{s_0}$, rather than $T_{w s_0} \cdot B_{s_0}$.
The step itself is straightforward: for each $x < w$, let
\[\TLLup(T_x, T_w) := p_{w, s_0} \circ \TLLup(T_w, T_{ws_0} \cdot B_{s_0}),\]
and for $x = w$ set
\[\TLLup(T_w, T_w) := p_{w, s_0} \circ \TLLup(T_{ws_0} \cdot B_{s_0}, T_{ws_0} \cdot B_{s_0}) \circ i_{w, s_0}.\]
Do similarly for the $\TLLdown$.

\subsection{Morphism pruning}

The previous step may have introduced some linear dependencies (over $R$) in the new $\TLLup$ and $\TLLdown$ sets.
These may be removed by standard means.

However, there is one way in which we can ``cheat'' a lot: since we are really only interested in holding on to these morphisms for the purposes of translating them and creating new local intersection forms, for $x \leq w$ we do not actually need the $R$-span of $\TLLup(T_x, T_w)$ to reflect the hom space $\Hom_{\cH_\bbO}(B_x, B_w)$, but rather we only need it to reflect the hom space $\Hom_{\cH_\bbO}^{\geq x}(B_x, B_w)$ (and similarly for the $\TLLdown$).
These are also precisely the sets we need in order to calculate the graded rank of $B_w$.
Therefore during this step we can in fact look at each element of $\TLLup(T_x, T_w)$ under the image of the linear map $\cat{Loc}^{\geq x}: \TLLup(T_x, T_w) \to \Hom_{\Std_Q}^{\geq x}(T_x, T_w)$, and extract a subset which descends to an $R$-independent set in the image.

Doing this of course changes the meaning of the sets $\TLLup$ and $\TLLdown$ over the whole algorithm, but provides a phenomenal speedup to the whole algorithm.

\subsection{Implementation details}

In order to implement this algorithm efficiently and with minimal fuss, we require several features commonly found in computer algebra systems (our prototype implementation is in MAGMA):
\begin{itemize}
    \item Implementations of polynomial rings over the integers or rationals, and their fraction fields (multivariate rational function fields).
    \item Sparse matrices with entries in arbitrary rings.
    \item A Coxeter group implementation capable of working with the Weyl and affine Weyl groups.
\end{itemize}

The \textit{translate morphisms} and \textit{braiding} steps of the algorithm rely on being able to act on objects and morphisms in $\Std_Q$ by objects and morphisms in $\cH_{BS}$, and the \textit{reduce idempotent} and \textit{morphism pruning} steps need to be able to project morphisms into the local quotient category $\Std_Q^{\not < w}$.
Our data structures for representing objects and elements of $\Std_Q$ closely matches the formal definition of an additive envelope:
\begin{itemize}
    \item
    A morphism $\varphi \colon X \to Y \in \Std_Q$ of degree $d$ is represented by a quadruple $(X, Y, [\varphi], d)$, where $X$ and $Y$ are sequences of elements in ${W}$, and $[\varphi]$ is the matrix of $\varphi$ as defined earlier, stored as a sparse matrix.
    We make this quadruple into its own datatype, with addition, scalar multiplication, and composition defined in the obvious way.
    For example, the composition $(X, Y, [\varphi], d_1) \circ (W, Z, [\psi], d_2)$ is $(W, Y, [\varphi][\psi], d_1 + d_2)$, where we additionally check that $X = Z$ and throw an error otherwise.
    Doing domain/codomain compatibility checks does not add much overhead to running time (at most a percent or two), and has the benefit of adding a lot of implict error checking through the program.
    \item
    A morphism $\psi \colon A \to B \in \cH_{BS}$ of degree $d$ is represented by a triple $(A, B, \Lambda(\psi))$, where $\Lambda(\psi) \in \Hom_{\cat{Std}}(\Lambda(A), \Lambda(B))$ is a standard morphism represented as described in the previous paragraph.
    If we intend to eventually act on ${^I \Std_{Q_I}}$ by $\psi$ we must remember the original Bott-Samelson domain and codomain $(A, B)$, it is not enough to remember only their images in $\cat{Std}$.
    \item
    The local quotient functor ${^I \Std_{Q_I}} \to {^I \Std_{Q_I}}^{\not <w}$ applied to the morphism $(X, Y, [\varphi], d)$ is $(X', Y', [\varphi]_{X', Y'}, d)$, where $X'$ is the subsequence of $X$ consisting of elements $\geq w$, and likewise for $Y'$, the matrix $[\varphi]_{X', Y'}$ is the submatrix defined by these subsequences.
 \end{itemize} 

With these data structures in place, it is fairly clear how to program the remainder of the algorithm.
We refer the reader to our implementation for further details.

\subsection{Limitations} \label{section:limitations}

The biggest limitation of this algorithm is the memory use: the fact that $\TLLdown(w, -)$ and $\TLLup(-, w)$ need to be stored somewhere in order to generate $\TLLdown(ws, -)$ and $\TLLup(-, ws)$ for $w < ws$.
For small primes especially, the structure of the Hecke category becomes very complex, and the sizes of these sets become infeasible to store, either in memory or on disk.
A limiting factor in particular is the sizes of coefficients\footnote{
    For this purpose, the \textit{size} of a rational number $a/b$ in reduced terms is $\log(|a| + |b|)$.
} encountered in the matrices, which seem to grow super-exponentially.
The basis of truncated light leaves is highly non-canonical, and the authors suspect that there may exist a basis which is better-suited to these calculations.


\section{A simpler algorithm} \label{section:simpler-algorithm}

Here we present a simpler algorithm for calculating the $p$-canonical basis, which works especially well for finite Coxeter groups of small rank.
It uses localisation to compute in the Hecke category, but in a much more straightforward way: rather than tracking the full set of light leaves associated to indecomposable objects, we merely use localisation to calculate local intersection forms of Bott-Samelson objects.
As the intersection forms can get extremely large, it is in our interest to calculate as few of them as possible, and to this end a whole host of known symmetries of $p$-canonical basis elements is taken advantage of.
The resulting algorithm is effective at calculating all $p$-canonical basis elements for groups up to rank 6 on a moderately powered desktop computer in less than a week, and up to rank 5 in minutes.

To fix notation, let $\set{\delta_w}$ denote the standard basis, $\set{b_w}$ the canonical basis, and $\set{{^p b}_w}$ the $p$-canonical basis of the Hecke algebra.
Define the coefficients $\set{h_{x, y}}$, $\set{{^p h}_{x, y}}$ and $\set{{^p m_{x, y}}}$ by the following changes of basis:
\begin{align}
    b_y &= \sum_{x \leq y} h_{x, y} \delta_y,
    &
    {^p b}_y &= \sum_{x \leq y} {^p h}_{x, y} \delta_x
    = \sum_{x \leq y} {^p m}_{x, y} b_x.
\end{align}
As the coefficients $h_{x, y}$ and ${^p h}_{x, y}$ are quite dense, in the sense that they are nonzero for all $x \leq y$, we will instead aim to instead calculate the coefficients ${^p m_{x, y}}$, or in other words we will work out the $p$-canonical basis relative to the canonical basis rather than the standard basis.

In \cref{section:local-intersection-form} we outlined the local intersection form: here we will instead only use the local intersection form inside the Bott-Samelson category $\cH_{BS}$.
The degree $d$ \defn{local intersection pairing} of $B_{\underline{x}}$ at $y \leq x$, after having fixed a reduced expression $\underline{y}$ for $y$, is the $R_\bbZ$-valued bilinear form
\begin{equation}
    \begin{aligned}
        I_{y, \underline{x}}^d
            \colon
            \Hom_{\not < y}^d(B_{\underline{x}}, B_{\underline{y}})
            \times
            \Hom_{\not < y}^{-d}(B_{\underline{y}}, B_{\underline{x}})
            &\to \End_{\not < y}^0(B_{\underline{y}}) \cong R_\bbZ, \\
        (f, g) &\mapsto f \circ g,
    \end{aligned}
\end{equation}
where the canonical isomorphism $\End_{\not < y}(B_{\underline{y}}) \to R_\bbZ$ is by taking the coefficient of the identity morphism in the double leaves basis for $\End(B_{\underline{y}})$.
The light leaves (not double leaves) $\LLdown(\underline{x}, \underline{y}) \subseteq \Hom_{\cH_{BS}}(B_{\underline{x}}, B_{\underline{y}})$ descend to a basis in the quotient $\Hom_{\not < y}(B_{\underline{x}}, B_{\underline{y}})$, and similarly for $\LLup(\underline{y}, \underline{x})$, so $I_{\underline{x}, y}$ may be viewed as a matrix whose rows and columns are in bijection with the subexpressions $\underline{e} \subseteq \underline{x}$ evaluating to $x^e = y$.

The \defn{graded rank} of the intersection pairing is
\begin{equation}
    \operatorname{\underline{\rank}} I^\bullet_{y, \underline{x}} = \sum_{d \in \bbZ} \left(\rank I^d_{y, \underline{x}} \right) v^d.
\end{equation}
Note that the graded rank is defined over $\bbZ$ because the Bott-Samelson category $\cH_{BS}$ is defined over $\bbZ$.

Now we introduce a $p$-modular system $(\bbO, \bbK, k)$: recall that the indecomposable object ${^\bbO B_x}$ base changes to the indecomposable object ${^\bbO B_x} \otimes k = {^k B_x}$ in the residue field, but in general the characteristic-zero object ${^\bbO B_x} \otimes \bbK$ is decomposable.
In terms of the bases above, we have $b_x = [{^\bbK B_x}]$ and ${^p b_x} = [{^k B_x}] = [{^\bbO B_x}]$, with the base change coefficients ${^p m_{x, y}} = [{^\bbK B_x}: {^\bbO B_y} \otimes \bbK]$.

The next theorem follows from \cite[Corollary~11.75]{eliasIntroductionSoergelBimodules2020} and idempotent lifting.

\begin{theorem}
    Let $(\bbO, \bbK, k)$ be a $p$-modular system.
    The multiplicity of the indecomposable object ${^\bbO B}_y$ in the Bott-Samelson object $B_{\underline{x}}$ is $[{^\bbO B}_y :B_{\underline{x}}] = \operatorname{\underline{\rank}} \left(I^\bullet_{y, \underline{x}} \otimes k\right)$.
\end{theorem}

With the technology developed for the previous algorithm we can certainly evaluate these intersection forms, and by the equations
\begin{align}
    [B_{\underline{w}}] &= [{^\bbO B_w}] + \sum_{x < w}[{^\bbO B_x} : B_{\underline{w}}][{^\bbO B_x}] \\
    b_{\underline{w}} &= {^p b_w} + \sum_{x < w} \operatorname{\underline{\rank}} \left(I^\bullet_{y, \underline{x}} \otimes k \right) b_x
\end{align}
we can calculate basis elements inductively\footnote{
    The Bott-Samelson element $b_{\underline{w}} = b_{w_1} \cdots b_{w_n}$ can be expanded back into the canonical basis quickly, provided one has a table of $\mu$-coefficients, i.e. the $W$-graph of the group.
}.
This forms the main idea of the algorithm, but we approach the computation slightly differently so that we can minimise the number of intersection forms calculated.

Split the coefficients up by degree: ${^p m}_{x, w} = \sum_{d \in \bbZ} {^p m}_{x, w}^d v^d$.
Whenever $ws < w$, then ${^\bbO B_w}$ is a summand of ${^\bbO B_{ws}} \otimes B_s$, and likewise for the left descents $sw < w$.
Let $U_{x, w}^d$ be the minimum coefficient of $b_x v^d$ appearing in the elements ${^p b_{ws}} \cdot b_s$ for $ws < w$ and $b_s \cdot {^p b_{sw}}$ for $sw < w$, then the $U_{x, w}^d$ give upper bounds on the $m_{x, w}^d$ coefficients.
Whenever $U_{x, w}^d = 0$ we have that ${^p m}_{x, w}^d = 0$: this becomes a \textit{known coefficient}.
There are a few other ways in which a coefficient might be known, such as the star operations mentioned below.

In order to calculate ${^p b_w}$, fix a reduced expression $\underline{w}$ and assume we have calculated ${^p b_x}$ for all $x < w$.
Fix an ordering $x_1, \ldots, x_r$ on the set of elements strictly less than $w$, such that $x_i \leq x_j$ in the Bruhat order implies $i < j$ (this order is necessarily weakly decreasing in length, and in fact any such order will do). We will calculate the sequence $\lambda_0, \ldots, \lambda_r$ defined by 
\begin{equation}
    \lambda_i = [B_{\underline{w}}] - \sum_{1 \leq j \leq i} [{^\bbO B_{x_i}} : B_{\underline{w}}] {^p b_{x_i}}, \quad \text{for } 0 \leq i \leq r.
\end{equation}
The first element $\lambda_0 = b_{\underline{w}}$ is the character of a Bott-Samelson object, and the last element $\lambda_r = [{^\bbO B_w}]$ is the character we are after.
When expanded in the canonical basis, the element $\lambda_i$ is equal to ${^p b_w}$ modulo lower terms in the canonical basis, specifically modulo $b_{x_j}$ for $j > i$ (this follows from the triangularity of the canonical basis and $p$-canonical basis).
Another way to say this is that the coefficient of $b_{x_j}$ in $\lambda_i$ is ${^p m}_{x_j, w}$ for $1 \leq j \leq i$.

At each step of the algorithm we need to determine the correct multiple $[{^\bbO B_{x_i}} : B_{\underline{w}}]$ of ${^p b_{x_i}}$ to subtract from $\lambda_{i-1}$ in order to get $\lambda_i$.
Let that multiple be $\mu_i$, then we need to determine the nonnegative integer coefficients $v^i$ inside $\mu_i$:
\begin{enumerate}
    \item If the coefficient of $v^d b_{x_i}$ inside $\lambda_{i-1}$ is zero, then $[v^d : \mu_i] = 0$.
    \item If ${^p m_{x_i, w}^d}$ is a known coefficient, then the coefficient of $v^d$ in $\mu_i$ must be ${^p m_{x_i, w}^d} - [b_{x_i}v^d : \lambda_{i-1}]$.
    \item Otherwise, calculate the rank of the local intersection form $I_{x_i, \underline{w}}^d \otimes k$, which is the multiplicity $[v^d \colon \mu_i]$ exactly.
\end{enumerate}

This strategy already gives a major improvement on the simple induction, due to the fact that far fewer intersection forms are calculated.
We use a few more symmetries of the $p$-canonical basis to speed it up further, which are outlined below.

By \cite[Proposition~4.2]{jensenCanonicalBasisHecke2016} we have ${^p m}_{x, y} = {^p m}_{x^{-1}, y^{-1}}$.
This means that once the element ${^p b_w}$ is calculated relative to the canonical basis, we can apply the inverse to every term to also get the element ${^p b_{w^{-1}}}$ for free whenever $w \neq w^{-1}$.
Similarly, if $\sigma \colon W \to W$ is a group automorphism induced by a Dynkin diagram automorphism, then ${^p m_{x, y}} = {^p m_{\sigma(x), \sigma(y)}}$.
So again we get an element for free whenever $w \neq \sigma(w)$, and furthermore the above symmetries can be combined so that yet another basis element comes for free whenever $w \neq \sigma(w^{-1})$.
Note that $\sigma$ commutes with inversion, and $\sigma$ is really only helpful in types $A$, $D$, and $E_6$.

Some extra symmetries come from the \textit{star operations}.
Given a pair $\set{s, t} \subseteq S$ with $m_{st} \in \{3, 4, 6\}$ there is a right star $w \mapsto w^*$ and left star $w \mapsto {^* w}$ involution, each defined on a subset of $W$.
By \cite[Corollary 4.7]{jensenABCPcells2020}, we have ${^p m_{x, y}} = {^p m_{x^*, y^*}}={^p m_{{^*x}, {^*y}}}$ whenever the prime $p$ is \textit{good} for $m_{st}$: order 3 is always good, order 4 is good for $p \geq 3$, and order 6 is good for $p \geq 6$.
This symmetry cannot be used to deduce ${^p b_{w^*}}$ from ${^p b_w}$ (since the involution is only partially defined), but it can be used to deduce a lot more \textit{known coefficients}, and thus reduce the number of intersection forms which need to be calculated. 

Finally, there are some shortcuts that can be taken when calculating the intersection form itself.
Rather than building a whole localised light leaf matrix, only to take its last column or row in the local $\not < x$ quotient, we can instead only calculate the last column or row.
The localised matrices are valued in the field of fractions of $\Sym(X)$, i.e. the multivariate fraction field $R(\alpha_1, \ldots, \alpha_n)$, but the intersection form will always turn out to be valued in $\bbZ$.
We can take another shortcut: we calculate the localised matrices over $\bbQ$, substituting each multivariate fraction $f$ for $f(1, 1, \ldots, 1)$.
Since a root is either positive or negative, evaluating it at $(1, \ldots, 1)$ will always yield a positive or negative number and never zero, therefore we get the right answer.

    \printbibliography

\end{document}